\RequirePackage{fix-cm}

\documentclass[twocolumn]{svjour3}          
  
\smartqed  

\usepackage{graphicx}
\usepackage{amsmath}
\usepackage{amssymb}
\usepackage[hidelinks]{hyperref}
\usepackage[noadjust]{cite}
\usepackage[ruled,vlined,linesnumbered]{algorithm2e}
\usepackage{xcolor}
\graphicspath{{Figures/}}

\begin{document}

\title{A Novel and Fully Automated Domain Transformation Scheme for Near Optimal Surrogate Construction}

\author{Johann Bouwer$^\dagger$  \and
        Daniel N. Wilke$^\dagger$$^*$ \and
        Schalk Kok$^\dagger$
}

\institute{   $\dagger$Department of Mechanical and Aeronautical Engineering, \\
              University of Pretoria, South Africa \\
              $*$\email{nico.wilke@up.ac.za}           
}

\date{Received: date / Accepted: date}

\maketitle

\begin{abstract}
Recent developments in surrogate construction predominantly focused on two strategies to improve surrogate accuracy. Firstly, component-wise domain scaling informed by cross-validation. Secondly, regression to construct response surfaces using additional information in the form of additional function-values sampled from multi-fidelity models and gradients. 

Component-wise domain scaling reliably improves the surrogate quality at low dimensions but has been shown to suffer from high computational costs for higher dimensional problems. The second strategy, adding gradients to train surrogates, typically results in regression surrogates. Counter-intuitively, these gradient-enhanced regression-based surrogates do not exhibit improved accuracy compared to surrogates only interpolating function values.

This study empirically establishes three main findings. Firstly, constructing the surrogate in poorly scaled domains is the predominant cause of deteriorating response surfaces when regressing with additional gradient information. Secondly, surrogate accuracy improves if the surrogates are constructed in a fully transformed domain, by scaling and rotating the original domain, not just simply scaling the domain. The domain transformation scheme should be based on the local curvature of the approximation surface and not its global curvature. Thirdly, the main benefit of gradient information is to efficiently determine the (near) optimal domain in which to construct the surrogate.

This study proposes a foundational transformation algorithm that performs near-optimal transformations for lower dimensional problems. The algorithm consistently outperforms cross-validation-based component-wise domain scaling for higher dimensional problems. A carefully selected test problem set that varies between 2 and 16-dimensional problems is used to clearly demonstrate the three main findings of this study.

\end{abstract}

\section{Introduction}

This paper develops and proposes a novel domain transformation scheme, completed as a preprocessing step, to improve the performance of surrogate models. The main application of surrogate models is in the field of Surrogate Based optimisation (SBO), where  computationally expensive simulations are replaced with a surrogate model to reduce the computational cost of optimisation.

Although many papers \cite{Vu2019, Koziel2011a, Cheng2020, Viana2021} make the claim that in this scenario gradient information of the function is not available, it is often not the case. Many papers \cite{Ryu1985,Olhoff1995a,sen,Parente2003} detail procedures to calculate the design sensitivities for functions that are computed using the Finite Element Method (FEM) or Computational Fluid Dynamics (CFD). Many finite element packages have adjoint sensitivities implemented, for example, Calculix \cite{Dhondt}. This gradient information can be calculated with respect to many different design variables to perform optimisation in a wide range of problems such as shape optimisation, thermodynamics, and vibration analyses \cite{Ryu1985, sen, KomKov1986, Balagangadhar2001, Newman1999}. Many SBO studies do make use of the available design sensitivities \cite{Laurenceau2008, Laurenceau2012, Laurent2019, Kampolis2004, Bouhlel2019, Koehler1996}, but often report either a small improvement in surrogate model accuracy, or consider the computational cost of these models prohibitive. 

Attempts to improve surrogate model accuracy include scaling the input or design domain. The current standard preprocessing strategy merely scales all dimensions of the input domain between 0 and 1 \cite{Vu2019, Viana2021}. This implicitly assumes that the underlying function is isotropic, i.e.\ the function behaviour is similar in all directions for an equivalent isotropic measure such as distance. In other words, the function is equally sensitive to variations in all input variables when evaluated at the same isotropic metric \cite{Lewis1998, Toal2008}. Attempts to discard the implicit isotropic assumption include
\begin{itemize}
    \item component-wise scaling of the domain, i.e.\ distinct scaling factors per dimension, as an attempt to recover isotropy after scaling \cite{Urquhart2020, Jones2001}, or
    \item adapting the surrogate model  to explicitly handle non-isotropic functions \cite{Viana2021, Bouhlel2019, Toal2008, Bouhlel2016}.
\end{itemize}
The problem with these strategies is that the designer is either left with an under-performing surrogate model, the implicit isotropic assumption \cite{Toal2008}, or the surrogate models become computationally intractable to construct for higher dimensional problems (typically $\geq$ 10) \cite{Bouhlel2019}. Therefore current efforts aim to decrease the computational cost of constructing models that attempt to embrace the non-isotropic nature of functions in typical engineering problems \cite{Bouhlel2019, Bouhlel2016, ChaeWilke}.

The main contribution of this paper is identifying and proposing a domain transformation scheme (scaling and rotation) that is an essential preprocessing step before the surrogate surface is constructed. It will be shown that the performance error of a surrogate model can be described as a summation of two distinct error sources. Firstly, there is the error associated with the sparsity of information, i.e.\ there are too few samples in the design space to fully capture the behaviour of the underlying function. The second, and the focus of this research, is the error associated with the mismatch between the domain the surrogate is constructed in and the implicit assumptions of standard surrogate models. It is shown in this research that even for simple low-dimensional problems the second domain-based error source can have a larger influence on the accuracy of a surrogate model than the more known or discussed sparsity error source. This second error can diminish to zero if the transformation scheme maps the original unsuitable design space to an ideal construction space for the surrogate model. Therefore, the optimal construction of a surrogate model is defined in this paper as one where the model is constructed in an ideal transformed domain.    

The proposed domain transformation scheme makes use of local Hessian estimates. Two options are proposed: estimate the Hessian from gradient information, or from function values. The transformation scheme based on gradient information is computationally efficient, as opposed to requiring a high-dimensional optimisation problem to be solved as with the component-wise Kriging scaling \cite{Jones2001}. Hence, scaling is not left to the user to identify but rather resolved automatically from the available data resulting in near optimal transformation. The benefits of this domain transformation (that includes scaling and rotation) include
\begin{itemize}
    \item gradient-enhanced surrogate models consistently outperform function value-only surrogate models,
    \item the improvement from function value-only surrogates to gradient-enhanced surrogates becomes more apparent for higher-dimensional problems,
    \item significantly fewer data points are required to construct similar quality surrogate models when no scaling or only basic scaling is performed.
\end{itemize} 
This paper is constructed as follows. Firstly, this paper offers a brief overview of SBO emphasising the basic steps in its implementation. This is followed by a detailed expos\'{e} of the construction and training of typical function-value based and gradient-enhanced surrogate models. From this breakdown the isotropic characteristic is discussed, as well as other lesser discussed characteristics, and a transformation scheme is developed. Lastly, the transformation scheme is assessed on test or benchmark problems, using both function and gradient-enhanced models of the same flexibility, before conclusions and recommendations are offered. 

\section{Surrogate Based optimisation}
\label{sec:SBO}
In general, the unconstrained optimisation problem attempts to find some vector of designs variables, $\boldsymbol{x}= [x_1, x_2, ..., x_n]^{\mbox{\tiny T}} \in \mathcal{R}^n$, that minimises some scalar function $F(\boldsymbol{x}):\mathcal{R}^n\rightarrow \mathcal{R}$. In many modern engineering optimisation problems, the evaluation of the function $F(\boldsymbol{x})$ often includes a computationally expensive simulation. A recent example of this is the work completed by Thapa and Missoum~\cite{Thapa2022}.

Therefore, the implementation of surrogate-based optimisation is often used when dealing with time-consuming simulations. A surrogate attempts to replace the expensive simulation with a computationally in-expensive model. The processes of surrogate-based optimisation is summarized into three phases \cite{Vu2019}: 
\begin{itemize}
    \item Phase 1: Select and evaluate a training set of design vectors.
    \item Phase 2: Use the training set to construct a surrogate model.
    \item Phase 3: Solve the surrogate optimisation problem and update the training set.
\end{itemize}
Phases 2 and 3 are then repeated until some termination condition is met. 

Phase 2, the construction of the surrogates, is the topic of interest in this paper. Therefore, an in-depth discussion of the processes used for this step in SBO is given in Section \ref{sec:Construction}.

Phases 1 and 3 are outside the scope of this paper. Hence, in this study, we merely sample points using the defacto-standard Latin Hyper-Cube sampling (LHS) without the space-filling condition enforced \cite{McKay1979}.

\section{Basic Surrogate Models}
\label{sec:Construction}
Surrogate models can be classified into function-value based, gradient-enhanced and gradient-only \cite{OptWilke}. Note that surrogate models that regresses through both function value and gradient information are referred to as either gradient-enhanced (GE) models \cite{Laurent2019, Bouhlel2019}, cooperative models (CO) \cite{Laurenceau2008, Laurenceau2012}, or first order (FO) models \cite{Koehler1996, OptWilke}. For the remainder of this research gradient-enhanced (GE) is used to describe surrogate models that regresses through both gradient and function value information.

Common function-value based surrogate models include Kriging Models, Radial Basis Functions (RBF) and polynomial surrogate models \cite{Viana2021, Vu2019, Urquhart2020, Toal2008, Bouhlel2016, Bouhlel2019}. Firstly, before gradient information can be included in these models, the more familiar function-value based models must be discussed and derived. Following these derivations, the gradient-enhanced models can be discussed.

\subsection{Function-Value Based Surrogate Models}

\subsubsection{Polynomial Surrogate Models}

The simplest surrogate model to implement is the polynomial model. An $n$-dimensional polynomial function of order $k$ can be expressed as
\begin{equation}
    f_{\mbox{\tiny poly}}(\boldsymbol{x}) =  \sum_{i=1}^n\sum_{j=1}^k W_{ij} \phi_{ij}(\boldsymbol{x}) + W_0,
    \label{eq:poly}
\end{equation}
where $W_{ij}$ is the weight associated with $j^{th}$ order of the polynomial in the $i^{th}$ dimension in the design space, and $\phi_{ij}(\boldsymbol{x})$ is the basis function
\begin{equation}
    \phi_{ij}(\boldsymbol{x}) = x_i^j.
\end{equation} 
The sampled point $\boldsymbol{x}$ is a column vector of size $n \time 1$. Therefore, there are $K = (k \times n)+1$ basis functions when we include the constant basis function.

Note that typically the coupling terms such as $x_ix_j$, $i\ne j$ are omitted from the basis functions. This omission is due to the exponential growth of the number of weights needed to fit the surrogate if these terms are included \cite{Viana2021, Vu2019}. 

If $p$ samples are taken from the design space then Equation~(\ref{eq:poly}) can be re-written as the system of equations
\begin{equation}
    \begin{bmatrix}
    f_1 \\
    f_2 \\
    \vdots \\
    f_p
    \end{bmatrix}_{p \times 1}
    =
    \begin{bmatrix}
    \phi_{kn}(\boldsymbol{x}_1) & \hdots & \phi_{11}(\boldsymbol{x}_1) & 1 \\
    \phi_{kn}(\boldsymbol{x}_2) & \hdots & \phi_{11}(\boldsymbol{x}_2) & 1 \\
    \vdots & \vdots & \vdots  & \vdots  \\
    \phi_{kn}(\boldsymbol{x}_p) & \hdots & \phi_{11}(\boldsymbol{x}_p) & 1 \\ \\
    \end{bmatrix}_{p \times K}
    \begin{bmatrix}
    W_{K} \\
    W_{K -1} \\
    \vdots \\
    W_{1}
    \end{bmatrix}_{K \times 1},
    \label{eq:PolySys}
\end{equation}
where $f_i$ is the function value at the sample location $\boldsymbol{x}_i$, and the weights ($W_{ij}$ and $W_0$) were assembled into a single vector and renumbered from 1 to $K$. This system of equations can then be expressed as
\begin{equation}
    \boldsymbol{f} = \boldsymbol{M}_p(\boldsymbol{x})\boldsymbol{W}_p.
    \label{eq:PolySysSimp}
\end{equation}
The training of a polynomial surrogate refers to the task of finding the optimal values of the weight vector $\boldsymbol{W}_p$. If $p = K$, $\boldsymbol{M}_p$ is square, and the weights can be solved using a linear algebra solver as long as $\boldsymbol{M}_p$ has full rank. More commonly $p > K$, therefore the system is over-determined \cite{OptWilke}, and needs to be computed using the least squares form
\begin{equation}
    {\boldsymbol{M}_p}^{\intercal}\boldsymbol{f} = {\boldsymbol{M}_p}^{\intercal} \boldsymbol{M}_p\boldsymbol{W}_p.
    \label{eq:LeastSquares}
\end{equation}

\subsubsection{Kriging Models}

Kriging, sometimes referred to as a Gaussian process, was first introduced by D.N Krige in 1951~\cite{Krige1951}. Kriging, unlike the other two surrogate models presented in this paper, models the underlying function using a  statistical approach. Jones~\cite{Jones2001} offers an intuitive derivation of the model which will be summarised and adapted here. 

Firstly, the model starts by assuming that the underlying function is a normally distributed process with a mean of $\mu$ and a variance of $\sigma$. Assuming the underlying function is continuous, the correlation between two values $f_i$ and $f_j$ at locations $\boldsymbol{x}_i$ and $\boldsymbol{x}_j$ in the $n$-dimensional design space is modelled mathematically with
\begin{equation}
    \mbox{Corr}[f_i, f_j] = \exp{\sum_{k = 1}^n \epsilon_k \biggl( (\boldsymbol{x}_i)_k - (\boldsymbol{x}_j)_k \biggr)^2}.
    \label{eq:KrigCorr}
\end{equation}
Here $(\boldsymbol{x}_i)_k$ refers to the $k$-th component of the location vector $\boldsymbol{x}_i$.

The distance measure, $((\boldsymbol{x}_i)_k- (\boldsymbol{x}_j)_k)^2$, intuitively implies that the correlation between two points will be high if points $\boldsymbol{x}_i$ and $\boldsymbol{x}_j$ are near one another in the design space, and decrease as the points are further from one another. The $\epsilon_k$ variable is a hyper-parameter that quantifies the dependency of the correlation on the distance in the $k$-th dimension in the design space. The covariance can then be expressed as 
\begin{equation}
    \mbox{Cov}(\boldsymbol{f}) = \sigma^2 \boldsymbol{M}_k,
    \label{eq:KrigCorrSimp}
\end{equation}
where $\boldsymbol{M}_k$ in an $p \times p$ matrix, where $p$ is the number of samples in the design space. Each component of $\boldsymbol{M}_k$ is given by Equation~(\ref{eq:KrigCorr}).

To estimate the ideal values of $\mu$, $\sigma^2$, and $\epsilon_k$ for all $n$ directions a vector of sampled values $\boldsymbol{f}$ of length $p$ is used in the equations
\begin{equation}
    \hat{\mu} = (\boldsymbol{I}\boldsymbol{M}_k^{-1}\boldsymbol{f})^{-1}{\boldsymbol{I}\boldsymbol{M}_k^{-1}\boldsymbol{I}},
    \label{eq:KrigMean}
\end{equation}
\begin{equation}
    \hat{\sigma^2} = \frac{1}{p}(\boldsymbol{f} - \boldsymbol{I}\hat{\mu})\boldsymbol{M}_k^{-1}(\boldsymbol{f} - \boldsymbol{I}\hat{\mu}),
    \label{eq:KrigVar}
\end{equation}
where $\boldsymbol{I}$ is the identity matrix. Equations~(\ref{eq:KrigMean}) and (\ref{eq:KrigVar}) depend on the matrix $\boldsymbol{M}_k$ which, in turn, depend of the vector of $\epsilon_k$ values. The optimal $\boldsymbol{\epsilon}$ is found by maximising the log-likelihood function
\begin{equation}
    L(\boldsymbol{\epsilon}) = - \frac{l}{2} \log(\sigma^2) - \frac{1}{2}\log(|\boldsymbol{M}_k|).
    \label{eq:LogLike}
\end{equation}
Strategies to solve the optimisation problem in Equation~(\ref{eq:LogLike}) is presented in Section \ref{sec:HyperParameter}.

The predicted value from the Kriging model of the underlying function at some location $\boldsymbol{x}^*$ is found with 
\begin{equation}
    f_{\mbox{\tiny Krig}}(\boldsymbol{x}^*) = \hat{\mu} + \boldsymbol{r}^{\intercal}\boldsymbol{M}_k^{-1}(\boldsymbol{f} - \boldsymbol{I}\hat{\mu}),
    \label{eq:Kriging}
\end{equation}
where $\boldsymbol{r}$ is a vector of the correlations computed by Equation~(\ref{eq:KrigCorr}), of the new point $\boldsymbol{x}^*$ and the previously sampled points:
\begin{equation}
    \boldsymbol{r} = 
    \begin{bmatrix}
    \mbox{Corr}(f^*, f_1) \\
    \mbox{Corr}(f^*, f_2) \\
    \vdots \\
    \mbox{Corr}(f^*, f_p) 
    \end{bmatrix}.
\end{equation}

\subsubsection{Radial Basis Function Surrogate Models}

Radial basis function surrogates refer to the family of surrogates that use a linear summation of basis functions that depend on a distance measure between two points. Popular options as basis functions include
\begin{itemize}
    \item Inverse quadratic: $\phi(\boldsymbol{x}, \boldsymbol{c}, \epsilon) = \frac{1}{1 + \epsilon||\boldsymbol{x} - \boldsymbol{c}||}$,
    \item Multi-quadratic: $\phi(\boldsymbol{x}, \boldsymbol{c}, \epsilon) = \sqrt{||\boldsymbol{x} - \boldsymbol{c}|| + \epsilon^2}$,
    \item Gaussian: $\phi(\boldsymbol{x}, \boldsymbol{c}, \epsilon) = e^{-\epsilon||\boldsymbol{x} - \boldsymbol{c}||^2}$,
\end{itemize}
where the variable $\epsilon$ is referred to as the shape parameter and the point $\boldsymbol{c}$ is the center of the basis function. The most widely used basis function is the Gaussian function \cite{Vu2019,Koziel2011a}. 
The RBF surrogate is expressed as a linear combination of $K$ basis functions
\begin{equation}
    f_{\mbox{\tiny RBF}} = \sum_{i=1}^K W_i\phi_i(\boldsymbol{x}, \boldsymbol{c}_i, \epsilon).
    \label{eq:BasisRBF}
\end{equation}
This equation becomes a system of equations similar to Equation~(\ref{eq:PolySysSimp})
\begin{equation}
   \boldsymbol{f} = \boldsymbol{M}_R(\boldsymbol{x}, \boldsymbol{c}, \epsilon)\boldsymbol{W}_R,
   \label{eq:RBFSysSimp}
\end{equation}
but now the matrix $\boldsymbol{M}_R$ becomes
\begin{equation}
\boldsymbol{M}_R = 
\begin{bmatrix}
    \phi(\boldsymbol{x}_1, \boldsymbol{c}_1, \epsilon) & \phi(\boldsymbol{x}_1, \boldsymbol{c}_2, \epsilon) & \hdots & \phi(\boldsymbol{x}_1, \boldsymbol{c}_K, \epsilon) \\
   \phi(\boldsymbol{x}_2, \boldsymbol{c}_1, \epsilon) & \phi(\boldsymbol{x}_2, \boldsymbol{c}_2, \epsilon) & \hdots & \phi(\boldsymbol{x}_2, \boldsymbol{c}_K, \epsilon) \\
    \vdots & \vdots & \vdots & \vdots \\
    \phi(\boldsymbol{x}_p, \boldsymbol{c}_1, \epsilon) & \phi(\boldsymbol{x}_p, \boldsymbol{c}_2, \epsilon) & \hdots & \phi(\boldsymbol{x}_p, \boldsymbol{c}_K, \epsilon) 
\end{bmatrix}.
\end{equation}
The remaining parameters of the surrogate include the number and locations of the centres $\boldsymbol{c}$ and the value of the shape parameter $\epsilon$. 

A popular choice for the centres is to select $p = K$, meaning that the number of centres is equal to the number of sampled points and to position the centres at the location of the sampled points. For this choice the matrix $\boldsymbol{M}_R$ becomes square and the weight vector can be solved directly from Equation~(\ref{eq:RBFSysSimp}). This is the method implemented for this research.

Some research implemented a fussy K-means clustering scheme to allocate the centres in the domain \cite{Vu2019, Koziel2011a}. From this scenario the system once again becomes over-determined and the least squares solution in Equation~(\ref{eq:LeastSquares}), where $\boldsymbol{M}_p$ and $\boldsymbol{W}_p$ now become $\boldsymbol{M}_R$ $\boldsymbol{W}_R$ respectively, is then used.

As with the Kriging hyper-parameter problem, the selection of a good shape parameter value $\epsilon$ for the RBF surrogate  is discussed in Section \ref{sec:HyperParameter}.

\subsection{Gradient-Enhanced Implementations}

Gradient-enhanced models are typically separated into two categories, namely direct and indirect approaches \cite{Laurent2019, Bouhlel2019, Bouhlel2016, OptWilke}. These approaches refer to the usage of the gradient information obtained at every sampled location.

\subsection{Indirect Gradient Enhancement}

The indirect gradient enhancement approach typically refers to some infill strategy. Using the 1st order Taylor series expansion,
\begin{equation}
    f(\boldsymbol{x}_i + \Delta x_k \boldsymbol{e}_k) = f(\boldsymbol{x}_i) + \frac{\partial f}{\partial x_k}\Delta x_k,
    \label{Eq:Infill}
\end{equation}
additional points are added to the dataset, where $\Delta x_k$ is the distance from the known sample point to the new infilled point and $\boldsymbol{e}_k$ is the unit vector in the direction along dimension $k$ \cite{Laurenceau2008, Chung}. The function value at the infilled point is not evaluated explicitly but estimated from Equation~(\ref{Eq:Infill}), by making use of the available gradient information. This strategy does not scale well with dimensionality, as for each point in the dataset, this method adds $n$ points. This means that the system can quickly become  ill-conditioned due to the closeness of the newly added points.

\subsection{Direct Gradient Enhancement}

The other approach to gradient-enhanced models is to directly include the gradients in the models themselves. This can either be done in an interpolating sense, such that the model directly interpolates both the function and gradient information at every point in the design space  \cite{Viana2021}, or in a regression sense, such that the model neither exactly fits the function or gradient information, but rather attempts to fit both in the least squares sense \cite{OptWilke}. 

A regression-based model is typically preferred to a fully interpolating model for two main reasons. Firstly, computational simulations that require discretisation and iterative solvers can result in noisy solutions. Therefore, if the model fits the solutions exactly the model may fit more to the noise in the data than to the underlying function. Secondly, a full interpolation matrix in either higher dimensional or densely sampled problems may become prohibitively large to solve, while a regression-based model can still offer useful results at a more reasonable computational cost. Therefore, regression-based derivations are offered in this section for the discussed surrogate models.

Another reason that regression models are preferred in this research is that the goal of the numerical investigations is to isolate the effect that the domain transformation has on the performance of the surrogate model. Therefore, the flexibility of the function and gradient-enhanced models are kept constant (by keeping the number and location of the centres the same), so that the only variable that is altered is the domain transformation strategy. The effect of increased flexibility in gradient-enhanced models, and how this increased flexibility is achieved, are outside the scope of this research. 

\subsubsection{GE Models}

Both the polynomial and the RBF surrogate models can be expanded to include gradient information in their construction. This can be done by first taking the gradient of their associated basis functions

\begin{equation}
    \frac{d \phi_{ij}(\boldsymbol{x})}{d \boldsymbol{x}} = j\boldsymbol{x}^{j - 1},
    \label{eq:DerPoly}
\end{equation}

\begin{equation}
    \frac{d\phi(\boldsymbol{x}, \boldsymbol{c}, \epsilon)}{d\boldsymbol{x}} = -2\epsilon \phi(\boldsymbol{x}, \boldsymbol{c}, \epsilon)(\boldsymbol{x} - \boldsymbol{c}),
    \label{eq:DerRBF}
\end{equation}
where Equations~(\ref{eq:DerPoly}) and (\ref{eq:DerRBF}) return column vectors of the gradients of the polynomial and RBF basis functions respectively.

A new system of equations can then be created from the gradient information at each sampled point for $p$ samples for the polynomial surrogate model

\begin{equation}
\begin{bmatrix}
    \frac{d f_1}{d \boldsymbol{x}} \\
    \frac{d f_2}{d \boldsymbol{x}} \\
    \vdots \\
    \frac{d f_p}{d \boldsymbol{x}}
    \end{bmatrix}
    =
    \begin{bmatrix}
    \frac{d\phi_{kn}(\boldsymbol{x}_1)}{d \boldsymbol{x}} & \hdots & \frac{d\phi_{11}(\boldsymbol{x}_1)}{d \boldsymbol{x}} & 0 \\
    \frac{d\phi_{kn}(\boldsymbol{x}_2)}{d \boldsymbol{x}} & \hdots & \frac{d\phi_{11}(\boldsymbol{x}_2)}{d \boldsymbol{x}} & 0 \\
    \vdots & \vdots & \vdots  & \vdots  \\
    \frac{d\phi_{kn}(\boldsymbol{x}_p)}{d \boldsymbol{x}} & \hdots & \frac{d\phi_{11}(\boldsymbol{x}_p)}{d \boldsymbol{x}} & 0 \\
    \end{bmatrix}
    \begin{bmatrix}
    W_K \\
    W_{K -1} \\
    \vdots \\
    W_{1}
    \end{bmatrix},
    \label{eq:DerPloyM}
\end{equation}
or the RBF surrogate model
\begin{equation}
    \begin{bmatrix}
    \frac{d f_1}{d \boldsymbol{x}} \\
    \frac{d f_2}{d \boldsymbol{x}} \\
    \vdots \\
    \frac{d f_p}{d \boldsymbol{x}}
    \end{bmatrix}
    =
    \begin{bmatrix}
    \frac{d \phi(\boldsymbol{x}_1, \boldsymbol{c}_1, \epsilon)}{d \boldsymbol{x}} & \frac{d\phi(\boldsymbol{x}_1, \boldsymbol{c}_2, \epsilon)}{d \boldsymbol{x}} & \hdots & \frac{d\phi(\boldsymbol{x}_1, \boldsymbol{c}_K, \epsilon)}{d \boldsymbol{x}} \\
    \frac{d \phi(\boldsymbol{x}_2, \boldsymbol{c}_1, \epsilon)}{d \boldsymbol{x}} & \frac{d\phi(\boldsymbol{x}_2, \boldsymbol{c}_2, \epsilon)}{d \boldsymbol{x}} & \hdots & \frac{d\phi(\boldsymbol{x}_2, \boldsymbol{c}_K, \epsilon)}{d \boldsymbol{x}} \\
    \vdots & \vdots & \vdots \\
    \frac{d \phi(\boldsymbol{x}_p, \boldsymbol{c}_1, \epsilon)}{d \boldsymbol{x}} & \frac{d\phi(\boldsymbol{x}_p, \boldsymbol{c}_2, \epsilon)}{d \boldsymbol{x}} & \hdots & \frac{d\phi(\boldsymbol{x}_p, \boldsymbol{c}_K, \epsilon)}{d \boldsymbol{x}} \\
    \end{bmatrix}
    \begin{bmatrix}
    W_K \\
    W_{K-1} \\
    \vdots \\
    W_{1}
    \end{bmatrix}.
    \label{eq:DerRBFM}
\end{equation}
These two systems, Equations~(\ref{eq:DerPloyM}) and (\ref{eq:DerRBFM}), can then be written as either
\begin{equation}
    \nabla \boldsymbol{f} = \boldsymbol{M}_{p-fo}\boldsymbol{W}_{p-fo},
\end{equation}
in polynomial model case, or
\begin{equation}
    \nabla \boldsymbol{f} = \boldsymbol{M}_{R-fo}\boldsymbol{W}_{R-fo},
\end{equation}
in the RBF model case. The subscript $-fo$ denotes that first-order information is used in the system. The gradient information can then be added to the original function-based systems, Equations~(\ref{eq:PolySysSimp}) and~(\ref{eq:RBFSysSimp}), to create a new system of equations 
\begin{equation}
    \begin{bmatrix}
    \boldsymbol{f} \\
    \nabla \boldsymbol{f}
    \end{bmatrix}
    =
    \begin{bmatrix}
    \boldsymbol{M}_{p} \\
    \boldsymbol{M}_{p-fo}
    \end{bmatrix}
    \boldsymbol{W}_{p-GE},
    \label{eq:PolyGE}
\end{equation}
in the polynomial case, or the in RBF case,
\begin{equation}
    \begin{bmatrix}
    \boldsymbol{f} \\
    \nabla \boldsymbol{f}
    \end{bmatrix}
    =
    \begin{bmatrix}
    \boldsymbol{M}_{R} \\
    \boldsymbol{M}_{R-fo}
    \end{bmatrix}
    \boldsymbol{W}_{R-GE}.
    \label{eq:RBFGE}
\end{equation}
The weight vector now contains the subscript $-GE$ to show that the weights solved from this system are for the gradient-enhanced versions of the surrogate models.

An important characteristic to note of the GE models is the size of the systems that need to be solved. In the function-value based models $p$ scalar samples are taken of the underlying function, creating a system of size $p \times K$, while in the GE models $p$ scalars and $p$ gradient vectors of size $n \times 1$ are sampled, creating a $(p + p \times n) \times K$ system. As the weight vector, $\boldsymbol{W}_{GE}$, is the same size, specifically $K \times 1$ in both the function and GE models, the models are of equal flexibility. The difference between the function and GE models is therefore that the GE models are constructed by regressing the model to the gradient information using the least squares formulation (similar to Equation~(\ref{eq:LeastSquares})). 

\subsubsection{GE-Kriging Models}

The derivation for direct gradient-enhanced Kriging is more complex than the derivations for the other two models discussed in this paper. Therefore, for the sake of brevity, and not to distract from the main contribution of this paper, the interested reader is referred to the literature  \cite{Laurenceau2008, Laurenceau2012, Bouhlel2019} for the complete mathematical description and implementation.

\subsection{Hyper-parameter Selection Strategies}
\label{sec:HyperParameter}

Unlike the polynomial surrogate model, the Kriging and RBF models require the optimisation or tuning of hyper-parameters. This optimisation sub-problem in SBO is a widely researched and discussed topic in literature. These two surrogate models can require vastly different algorithms to find the optimal hyper-parameter or set of hyper-parameters. 

It has been shown that the numerical value of the hyper-parameter greatly impacts the performance of the model. Therefore, before further research can be completed it is necessary to discuss the current optimisation methods implemented for the hyper-parameter selection of these models.

\subsubsection{The Kriging Hyper-parameter problem}
\label{sec:KrigHyp}

The main challenge when solving the Kriging hyper-parameter optimisation problem in Equation~(\ref{eq:LogLike}), is the fact that there are as many $\epsilon$ values as there are dimensions in the sampled design space. 

Therefore many papers apply some global optimiser to solve this problem, such as the Genetic Algorithm (GA) or Particle Swarm Optimisation (PSO) \cite{Toal2008}. In higher dimensions, this becomes computationally expensive, so much so that it can become the bottleneck in computation time for SBO. Toal et al.~\cite{Toal2008} investigated four different tuning strategies on problems varying from 1D to 30D. Each of the tuning strategies sampled the model 10 000 times before a set of hyper-parameters was selected. 

Other papers attempt to reduce the number of hyper-parameters in the model. Bouhel et al.\ \cite{Bouhlel2019, Bouhlel2016} used a partial-least squares (PLS) method to introduce new kernels based on the information from the PLS method. The number of hyper-parameters is then reduced to the number of principal components (PC) the designer decides to keep based on the information gathered from the PLS method. The ideal number of PC to be retained depends on the problem as well as the location of the sampled points. There is currently no consistent method to determine this value.

The last option is to reduce the hyper-parameter vector to one value, i.e.\ one constant value for all the directions. This has been shown  \cite{Toal2008} to produce better results than the other two methods if the underlying function is isotropic in nature. This isotropic assumption significantly affects the accuracy of surrogates and will be discussed further in Section~\ref{sec:Charateristics}.

In this research a simplex search algorithm, such as that used by Toal et al. \ \cite{Toal2008}, is implemented to find optimum scaling values for the Kriging hyper-parameter problem. To keep the computationally costs reasonable, as well as competitive with the other methods implemented, the algorithm is limited to 100 iterations for 5 initial scaling vectors. 

\subsubsection{The RBF Hyper-parameter problem}

As in the case of Kriging, the selection of the shape parameter of the RBF model is an often discussed and researched topic in literature. Some papers propose some heuristic to calculate the single scalar value $\epsilon$, typically based on the dimensionality of the problem and the distance between the sampled points \cite{Golbabai2015,Benoudjit2003}. Others implement some cross-validation schemes such as K-fold cross-validation or leave-out-one cross-validation (LOOCV) \cite{Urquhart2020, Viana2021}.

In this research, LOOCV is implemented to optimise the shape parameter of the RBF model. The LOOCV method typically involves the following steps:
\begin{enumerate}
    \item Divide the dataset into many subsets, where each subset contains all the points except one (a different one for each subset).
    \item Construct the surrogate for each subset of the sampled points.
    \item Find the error for each constructed surrogate at the point that was excluded from the subset.
    \item Sum all the errors.
\end{enumerate}
These steps are then repeated for each trail shape parameter value in some predetermined set of values. The value that results in the lowest summed error value is then used to construct the surrogate on the full dataset. For large datasets, this can become computationally expensive as $p$ surrogates are trained for each tested shape parameter value. Therefore, the algorithm proposed by Rippa \cite{Rippa1999} is used. An estimated error value $E$ is computed from
\begin{equation}
    E = \frac{1}{K}\sum_{i = 1}^{K}\frac{{W}_i}{{M}^{-1}_{ii}},
    \label{eq:LOOCV}
\end{equation}
where the values of ${W}_i$ and ${M}^{-1}_{ii}$ are found by constructing the surrogate on the entire data set. $W_i$ is the solved weight for the $i$-th basis function, and $M_{ii}^{-1}$ is the $i$-th diagonal of the inverse of the basis function matrix. Equation~(\ref{eq:LOOCV}) allows for only one surrogate to be constructed per trail shape parameter value instead of $p$ surrogates. 

\section{Implicit Isotropic Assumption}
\label{sec:Charateristics}

\begin{figure}[h!]
    \centering
    \includegraphics[width = \columnwidth]{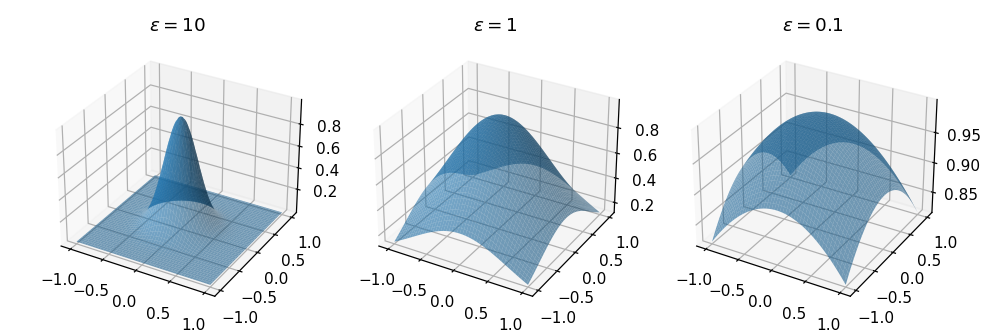}
    \caption{The Gaussian basis function for shape parameters of 10, 1, and 0.1.}
    \label{fig:Guass}
\end{figure}

This section discusses what is meant by isotropic and non-isotropic functions as well as why this characteristic can be detrimental to the performance of the surrogate. Section \ref{sec:Construction} demonstrated that surrogates are a linear combination of basis functions. Clearly, if these basis function shapes do not share some similarity with the underlying function, the performance of the overall surrogate will suffer. Figure~\ref{fig:Guass} illustrates the Gaussian basis function, the most common basis function, with three different shape parameter values for the two-dimensional case.

From the sub-figures in Figure \ref{fig:Guass} it is clear that the Gaussian basis function is symmetrical or isotropic. This means that the surrogate makes the implicit assumption that the variables of the underlying function are all equally important to the outcome of the function. 

This can be investigated in more depth by looking at the effect that the shape parameter has on the Gaussian basis function. Notice from Figure~\ref{fig:Guass} that changing the shape parameter $\epsilon$ only impacts the curvature of the basis function. Both the function value and the gradient vector are independent of the shape parameter at the center $\boldsymbol{x} = \boldsymbol{c}$. This can be seen mathematically in Equations~(\ref{eq:BasisRBF}) and (\ref{eq:DerRBF}) if the the equations are evaluated at the point $\boldsymbol{x} = \boldsymbol{c}$
\begin{equation}
    \phi(\boldsymbol{x}, \boldsymbol{c}, \epsilon) \vert_{\boldsymbol{x} = \boldsymbol{c}} = 1,
\end{equation}
\begin{equation}
    \biggl. \frac{d\phi(\boldsymbol{x}, \boldsymbol{c}, \epsilon)}{d\boldsymbol{x}}\biggr\vert_{\boldsymbol{x} = \boldsymbol{c}} = \boldsymbol{0}.
\end{equation}

The second derivative of the function is given by
\begin{equation}
 \frac{d^2\phi(\boldsymbol{x}, \boldsymbol{c}, \epsilon)}{d\boldsymbol{x}^2} = -2\epsilon\frac{d\phi}{d\boldsymbol{x}}(\boldsymbol{x} - \boldsymbol{c})^{\tiny{T}} -2\epsilon \boldsymbol{I} \phi(\boldsymbol{x}).
    \label{eq:HessBasis}
\end{equation}

If the second derivative is evaluated at the point $\boldsymbol{x} = \boldsymbol{c}$, this results in
\begin{equation}
   \biggl. \frac{d^2\phi(\boldsymbol{x}, \boldsymbol{c}, \epsilon)}{d\boldsymbol{x}^2}\biggr\vert_{\boldsymbol{x}=\boldsymbol{c}} = -2\epsilon\boldsymbol{I}.
    \label{eq:HessBasisSimp}
\end{equation}
Therefore, the act of altering or optimising the shape parameter therefore clearly results in an equal change in the curvature of the basis function (at the center) in all directions. 

\begin{figure}[h!]
    \centering
    \includegraphics[width = \columnwidth]{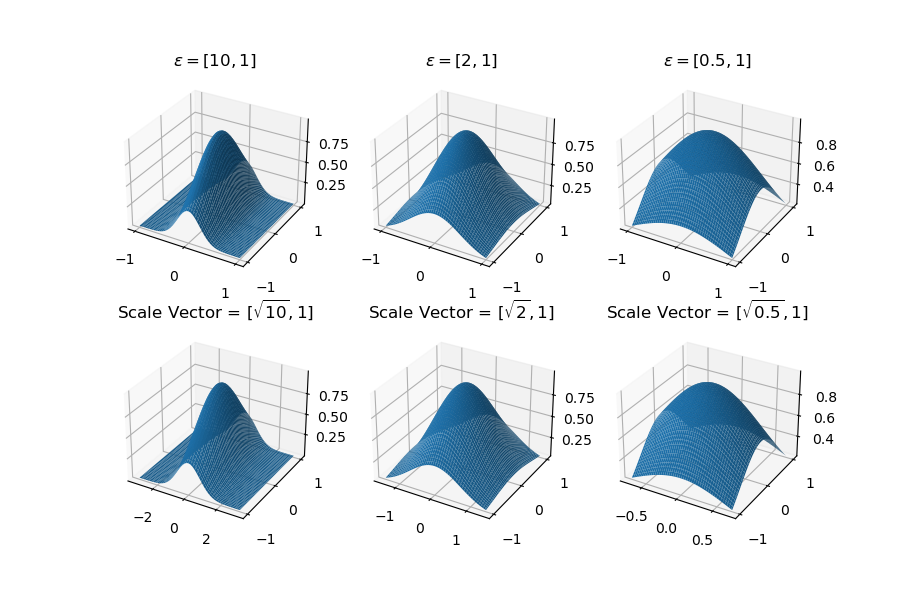}
    \caption{The Gaussian basis function with independent shape or scale parameters for each dimension.}
    \label{fig:InGauss}
\end{figure}

Therefore, if one shape parameter is used for all directions, the model will have ideal performance if the underlying function exhibits similar curvature in all directions. However, it is unlikely that a practical engineering design or optimisation problem will contain variables that all have equal (or at least similar) impact on the outcome of the design. Therefore, what is currently done is either a different shape parameter is assigned to each principal direction in the design space or a different scale parameter is used for each principal direction in the design space, such that the underlying function becomes isotropic. Figure~\ref{fig:InGauss} shows the Gaussian basis function with different shape parameters and scale parameters in each principal direction.

From Figure \ref{fig:InGauss} it can be seen that the two options are equivalent, i.e.\ for a certain shape parameter there is a corresponding scale parameter, specifically there is a square root relationship between the two, that will alter the curvature by either ``stretching'' or ``shrinking'' the domain such that the curvature will be equivalent.  This vector of hyper-parameters, either shape or scale parameters, means that the curvature in each principal direction can be altered independently, thus removing the implicit assumption that all variables impact the outcome equally.

As has already been discussed in Section~\ref{sec:Construction}, finding the optimal values of either scale or shape parameters for each dimension in a design problem, creates a computationally expensive hyper-parameter optimisation problem. However, what has not been discussed is that this formulation also makes the implicit assumption that the variables all independently impact the underlying function. This assumption is revealed from the observation that the shape parameters are incapable of changing the curvature in any direction other than the principal directions. This feature makes the implicit assumption that the variables all \emph{independently} influence the outcome of the function and that there is no inherent or underlying relationship between two or more variables. This can once again be seen by taking Equation~(\ref{eq:HessBasis}) and adapting it to an $n$-dimensional vector of $\epsilon$ values. The second derivative of the basis function is then given by
\begin{equation}
    \frac{d^2\phi(\boldsymbol{x}, \boldsymbol{c}, \epsilon)}{d x_n^2} = -2\epsilon_n.
    \label{eq:HessKrig}
\end{equation}
Notice that the Hessian at $\boldsymbol{x} = \boldsymbol{c}$ can be written as a diagonal $n \times n$ Hessian where the diagonal vector is $-2 \boldsymbol{\epsilon}$, i.e.\ a constant times the shape parameter vector. 

Therefore, the goal of this paper is to develop a transformation scheme that will define a domain that will force the underlying function to be \emph{both} isotropic and independent (uncoupled or decomposable). The surrogate can then be constructed in this transformed domain. The transformation scheme can be used to map to and from the original and transformed domains. This transformation scheme needs to be general and robust, 
i.e.\ it cannot return a domain in which the accuracy of the model will worsen, and it must be computationally efficient.

\subsection{Example Function}
\label{sec:NumericalExample}

To demonstrate the above arguments, the following uncoupled 2D function with each dimension in the domain $x_i \in [0,1]$ is considered:
\begin{equation}
    F(\boldsymbol{x}) = \sin(2\pi x_1) + \sin(2\pi x_2).
    \label{eq:NumericalExample}
\end{equation}
The effect that the scaling and rotating of the domain have on the performance of the RBF surrogate model is demonstrated by defining two new domains. Firstly, a scaled domain $\boldsymbol{x^*}$ is defined in which the domain of the function is scaled using the equation
\begin{equation}
    \boldsymbol{x^*} = \boldsymbol{S}\boldsymbol{x},
\end{equation}
where the matrix $\boldsymbol{S}$ is defined as 
\begin{equation}
    \boldsymbol{S} = 
    \begin{bmatrix}
    2 & 0 \\
    0 & 1
    \end{bmatrix}.
\end{equation}
The scaled domain $\boldsymbol{x^*}$ is then rotated to the domain $\hat{\boldsymbol{x}}$. In this domain, the function becomes coupled. The domain transformation is given by
\begin{equation}
    \boldsymbol{\hat{\boldsymbol{x}}} = \boldsymbol{R}\boldsymbol{x^*} = \boldsymbol{R} \boldsymbol{S} \boldsymbol{x}
\end{equation}
where the rotation matrix $\boldsymbol{R}$ is defined as 
\begin{equation}
    \boldsymbol{R}
    =
    \begin{bmatrix}
    \cos(30^\circ) & -\sin(30^\circ) \\
    \sin(30^\circ) & \; \; \; \cos(30^\circ)
    \end{bmatrix}.
\end{equation}

\begin{figure}[h!]
    \centering
    \includegraphics[width = \columnwidth]{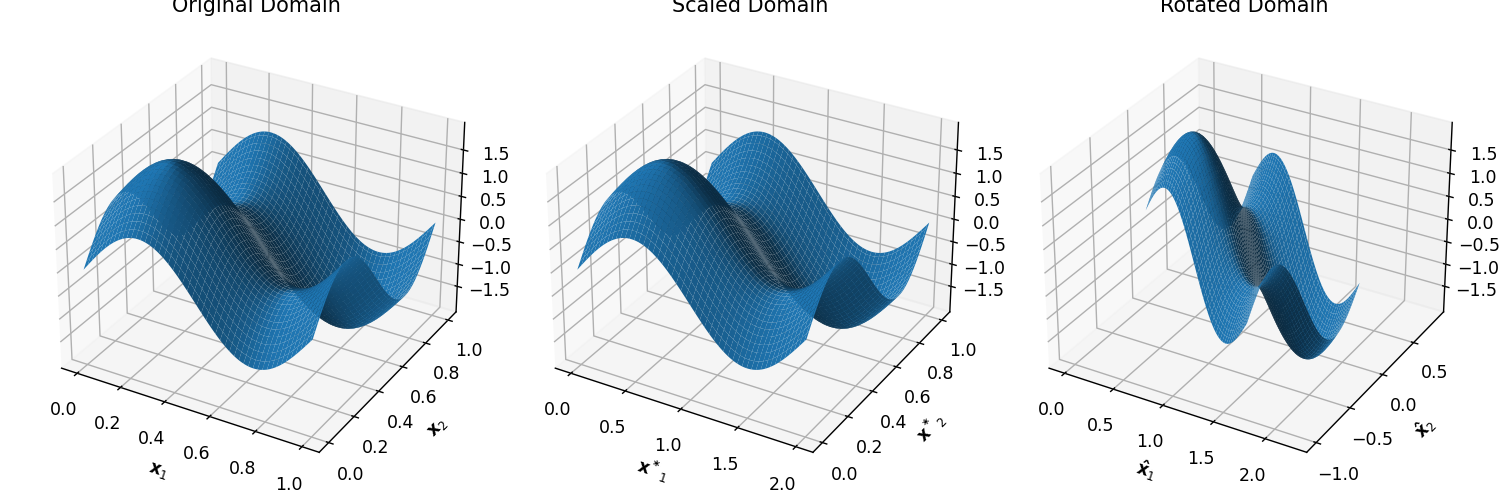}
    \caption{Example function illustrated in the original domain, scaled domain and rotated domain.}
    \label{fig:NumericalExample}
\end{figure}

\begin{figure}[h!]
    \centering
    \includegraphics[width = \columnwidth]{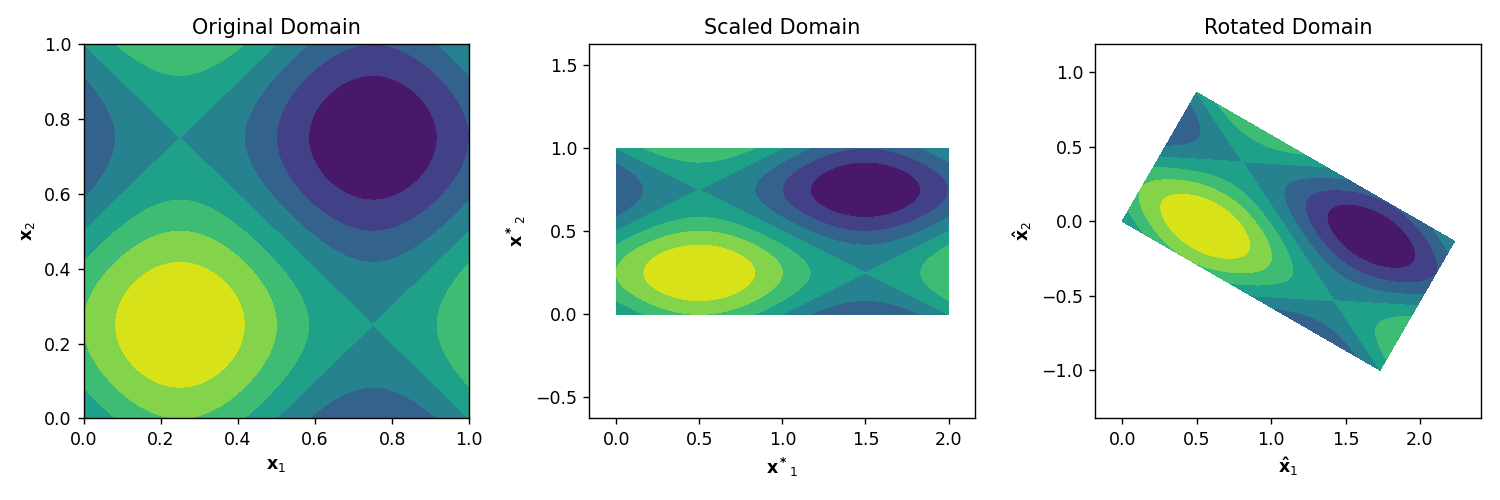}
    \caption{Contour plots of the example function illustrated in the original domain, scaled domain and rotated domain.}
    \label{fig:NumericalExampleContour}
\end{figure}

The function in the three domains, namely the original, scaled, and rotated domains is shown in Figures \ref{fig:NumericalExample} and \ref{fig:NumericalExampleContour}.

Three RBF surrogates are then constructed using various sample numbers (varying from 10 to 26), one in the original domain, one in the scaled domain, and lastly one in the rotated domain.

The performance of each surrogate is measured at 1000 randomly sampled test points. The number of test points is selected so much higher than the number of construction points to ensure that the error measure is an accurate reflection of the quality of fit, and is not affected by the location of the test points. To account for the randomness present in the location of the construction points, the error calculation is repeated 50 times and the mean is recorded. To evaluate the dependency of the surrogates on the locations of the construction points, a measure of the variance of the shapes of the surrogates is recorded. This is done by taking the variance of the error for each point in the test set and then recording the mean of this variance across all the points. Ideally, this result should be zero, otherwise, the surrogate greatly depends on the randomness of the sampling technique.

\begin{figure}[h!]
    \centering
    \includegraphics[width = \columnwidth]{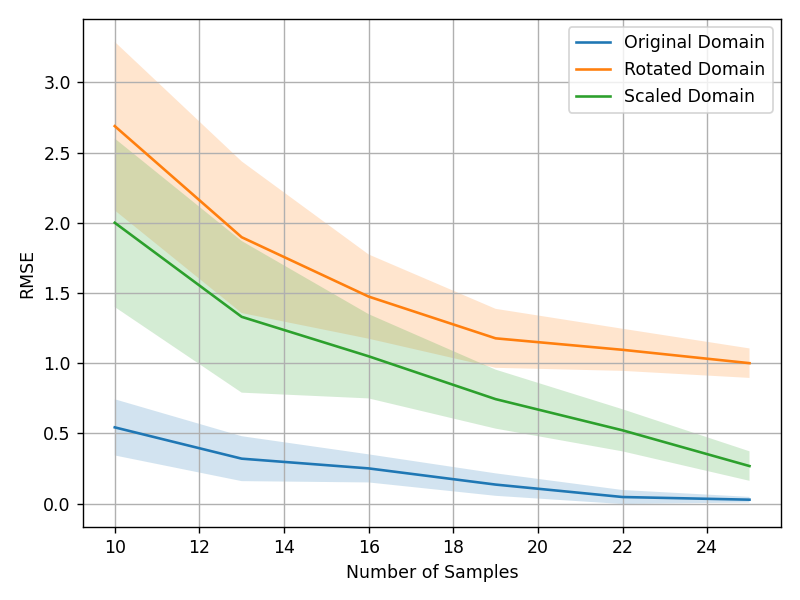}
    \caption{The mean (solid lines) and variance (shaded regions) in the RMSE for the surrogates constructed in the three domains for an increasing number of samples.}
    \label{fig:ExampleResFull}
\end{figure}

The performance measure used is the Root Mean Square Error (RMSE), expressed by
\begin{equation}
    RMSE = \sqrt{\frac{\sum_i^N(V_T^i - V_P^i)^2}{N}}
    \label{eq:RMSE}
\end{equation}

where $V_T^i$ is the target value and $V_P^i$ is the predicted value from the surrogate. The results are shown in Figure \ref{fig:ExampleResFull} where the shaded region in the plots indicates the variance of the surrogates.

Clearly, the domain the surrogate is constructed in has a meaningful and measurable impact on the performance of a surrogate. The transformed domains, i.e $\boldsymbol{x^*}$ and $\hat{\boldsymbol{x}}$, negatively impacted both the performance of the surrogate (increased error), as well as the consistency of the surrogate (increased variance), especially at lower sampling densities. One can also see the benefit of a complete transformation (rotation and scaling) that would transform the problem back from the rotated $\hat{\boldsymbol{x}}$ domain to the original $\boldsymbol{x}$ domain.

\begin{figure}[h!]
    \centering
    \includegraphics[width = \columnwidth]{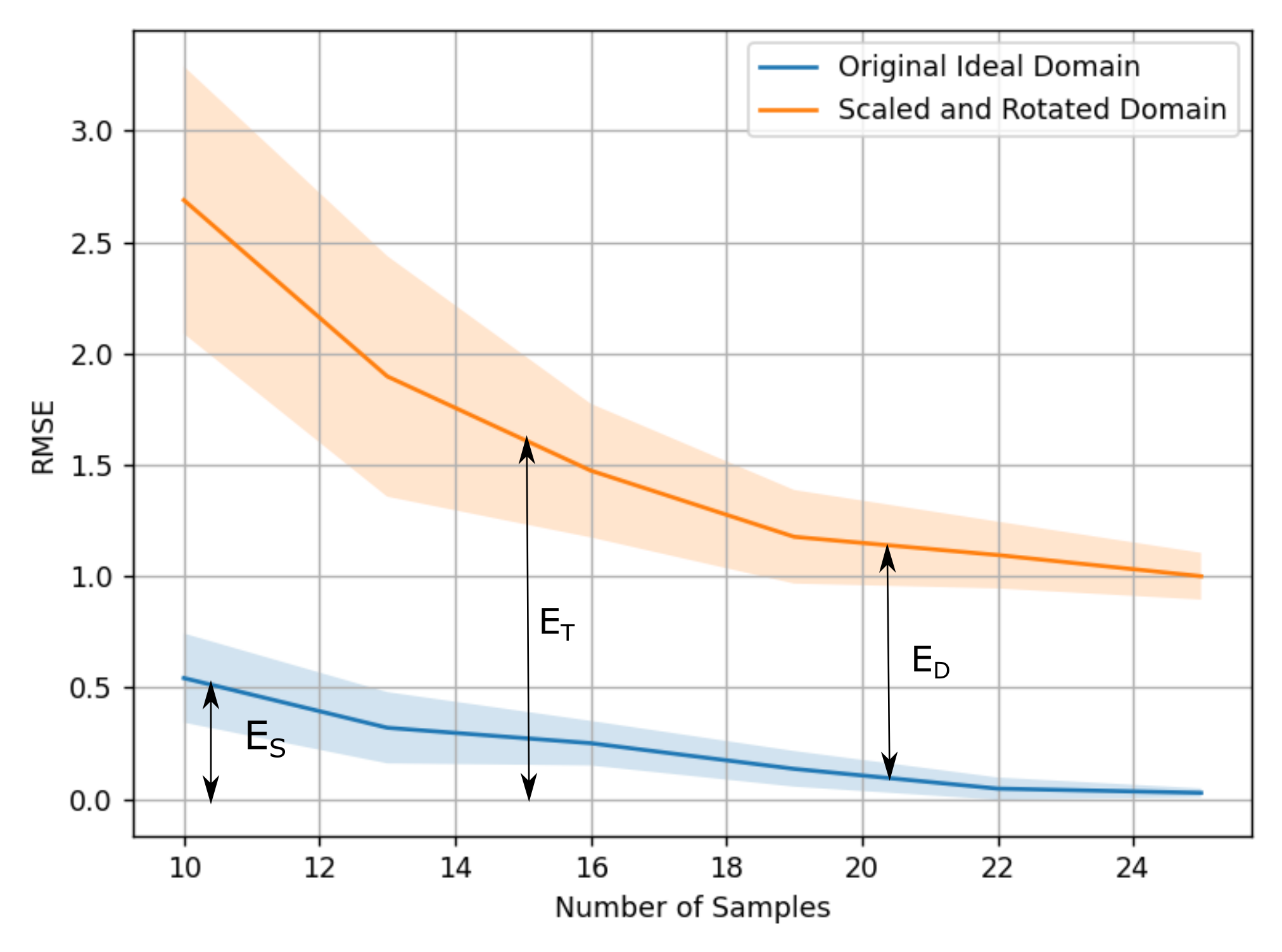}
    \caption{The sources of poor performance of a surrogate. The total error $E_T$ consists of the sparsity of information error $E_S$ and the construction domain error $E_D$. }
    \label{fig:ExampleRes}
\end{figure}

The total error of a surrogate, $E_T$, can then be defined as a summation of two errors. The first is the error associated with the sparsity of information, $E_S$, and the second is the error associated with the domain the surrogate is constructed in, $E_D$. These errors are indicated in Figure~\ref{fig:ExampleRes}.

\section{Proposed Transformation Scheme}
\label{sec:TransScheme}

In practical design or optimisation problems the ideal values of the matrices $\boldsymbol{R}$ and $\boldsymbol{S}$, the rotation and scaling matrices, will be unknown. Therefore, in this section, an efficient, consistent, and general domain transformation scheme is developed. To begin this discussion consider a simple multidimensional non-linear polynomial function, the 2-dimensional quadratic function. This function can be expressed as
\begin{equation}
    f(\boldsymbol{x}) = \frac{1}{2}\boldsymbol{x}^\intercal \boldsymbol{A}\boldsymbol{x} + \boldsymbol{b}^\intercal\boldsymbol{x} + c,
    \label{eq:quad}
\end{equation}
where $\boldsymbol{A}$ is a $2\times2$ matrix, $\boldsymbol{b}$ is a $2\times1$ vector, and $c$ is a scalar. For this discussion, the case where $\boldsymbol{b}$ and $c$ are zero is considered. Therefore Equation~(\ref{eq:quad}) becomes
\begin{equation}
    f(\boldsymbol{x}) = \frac{1}{2}\boldsymbol{x}^\intercal \boldsymbol{A} \boldsymbol{x}.
\end{equation}
In this form, the $\boldsymbol{A}$ matrix is equal to the Hessian or curvature matrix of the function.

\begin{figure}[h!]
    \centering
    \includegraphics[width = \columnwidth]{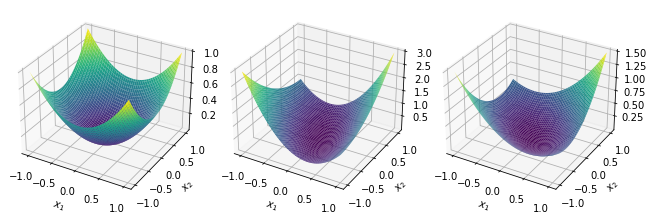}
    \caption{3D plots of quadratic $A_1$, $A_2$, and $A_3$ functions respectively.}
    \label{fig:Quad3D}
\end{figure}

\begin{figure}[h!]
    \centering
    \includegraphics[width = \columnwidth]{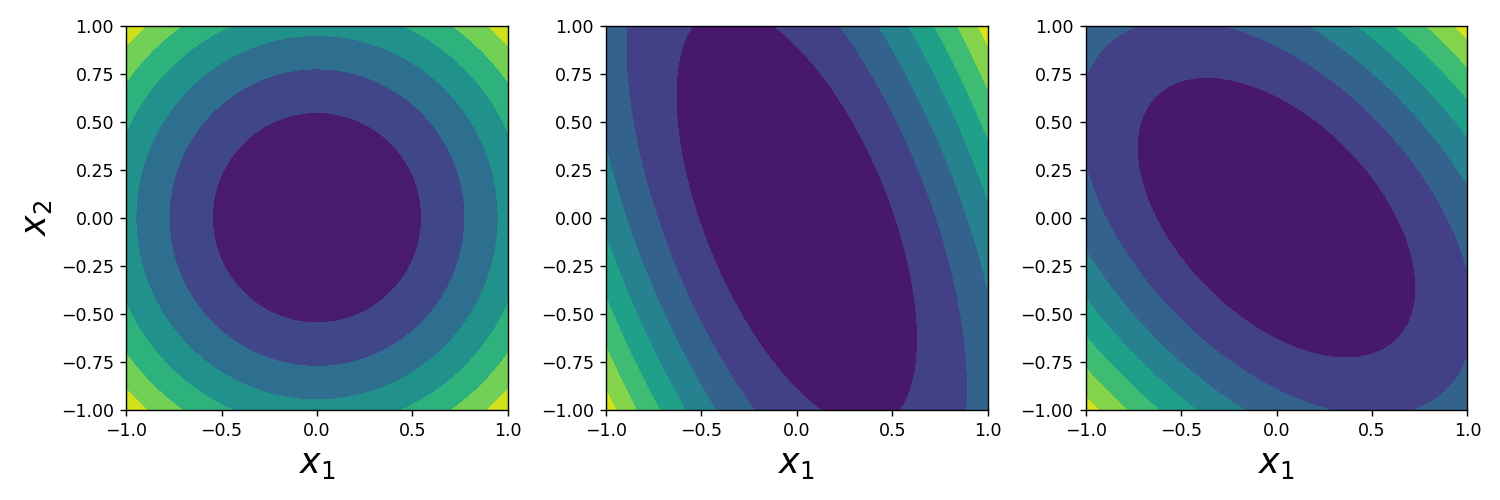}
    \caption{Contour plots of quadratic $A_1$, $A_2$, and $A_3$ functions respectively.}
    \label{fig:QuadContour}
\end{figure}

Figure~\ref{fig:Quad3D} shows the 3D representation and Figure~\ref{fig:QuadContour} shows the contour plots, for $f$, for three cases
\begin{equation}
    \boldsymbol{A}_1 = 
     \begin{bmatrix}
        1 & 0 \\
        0 & 1
    \end{bmatrix}, \; \;
    \boldsymbol{A}_2 = 
    \begin{bmatrix}
        3 & 1 \\
        1 & 1
    \end{bmatrix}, \; \;
    \boldsymbol{A}_3 = 
    \begin{bmatrix}
        1 & 0.5 \\
        0.5 & 1
    \end{bmatrix}.
\end{equation}

The shape of the function in the case of $\boldsymbol{A}_1$, when the function is isotropic, closely resembles the shape of the Gaussian basis function. Therefore, the Gaussian basis function is more suitable for the case where $\boldsymbol{A} = \boldsymbol{A}_1$ than when $\boldsymbol{A} = \boldsymbol{A}_2$ or $\boldsymbol{A} = \boldsymbol{A}_3$.

\begin{figure}[h!]
    \centering
    \includegraphics[width = \columnwidth]{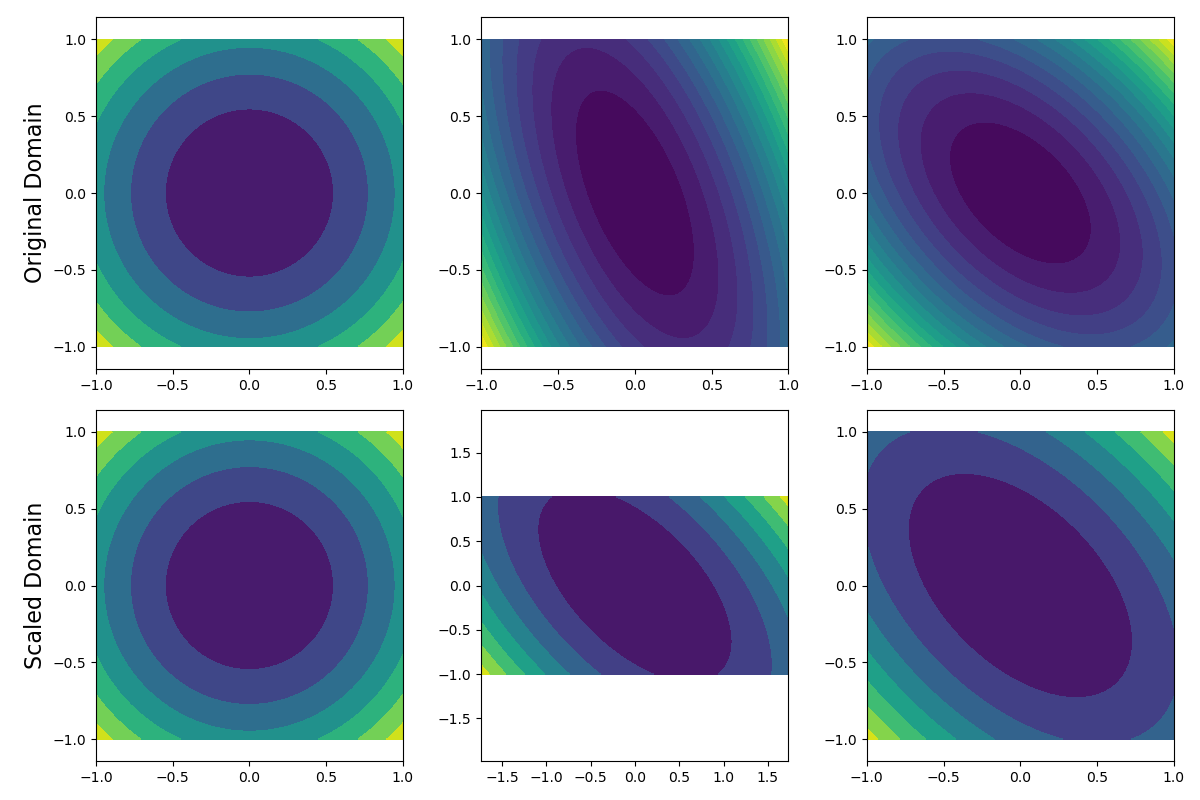}
    \caption{Quadratic functions scaled in the principal directions with the square root of the diagonal entries of their Hessians.}
    \label{fig:ScaledQuad}
\end{figure}

Therefore, the goal of the transformation scheme should be to create a domain where for any $\boldsymbol{A}$, the function evaluated in the transformed domain should resemble the case where $\boldsymbol{A} = \boldsymbol{A}_1$. If, as is currently a popular choice, the domain is only scaled independently in each principal direction, and there is coupling between variables (i.e.\ the Hessian matrix is not a diagonal matrix), then Figure~\ref{fig:ScaledQuad} is obtained. Here each dimension is scaled by the square root of the corresponding diagonal entry in the Hessian. 

Figure~\ref{fig:ScaledQuad} clearly demonstrates that only co-ordinate based scaling is insufficient to create an isotropic function. The relationship between the variables must therefore be taken into account.

\begin{figure}[h!]
    \centering
    \includegraphics[width = \columnwidth]{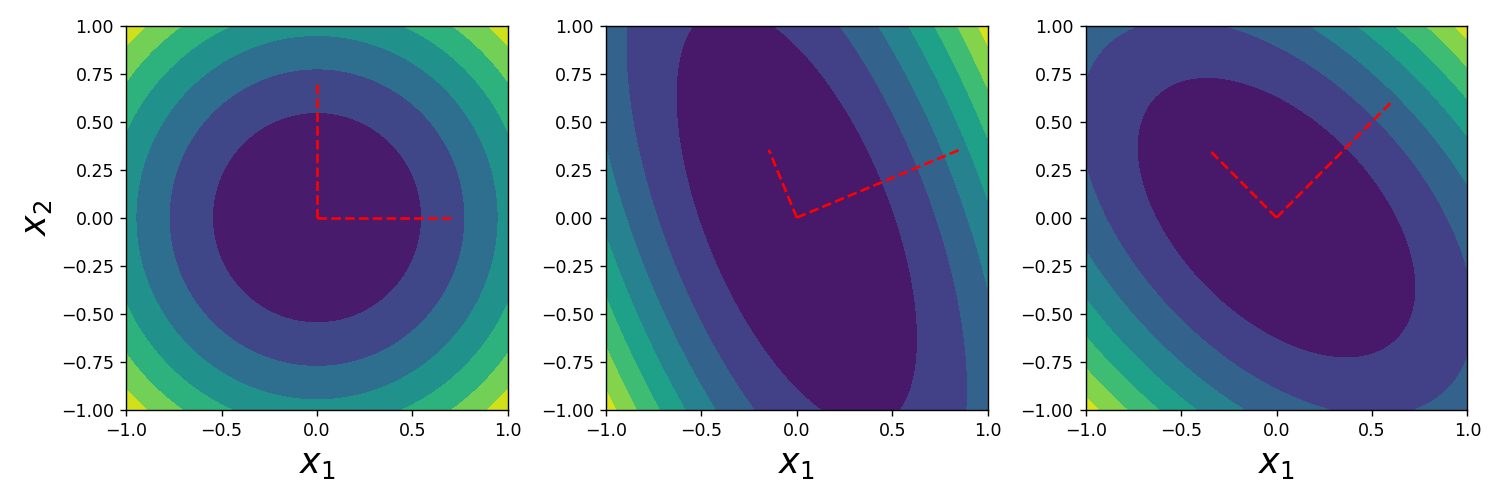}
    \caption{Contour plots of quadratic functions with Hessians given by $\boldsymbol{A}_1$, $\boldsymbol{A}_2$ and $\boldsymbol{A}_3$. The dashed lines indicate the eigenvectors and their lengths are chosen proportional to the eigenvalues.}
    \label{fig:QuadEig}
\end{figure}

A transformation scheme that can transform the domain such that the resulting function becomes isotropic can be achieved by considering the eigenvalues and eigenvectors of the Hessian. Figure~\ref{fig:QuadEig} shows the eigenvalues and eigenvectors for the three problems overlaid with their corresponding contour plots. The eigenvectors are indicated with dashed lines, with the length of each dashed line chosen proportional to the magnitude of the corresponding eigenvalue.

\begin{figure}[h!]
    \centering
    \includegraphics[width = \columnwidth]{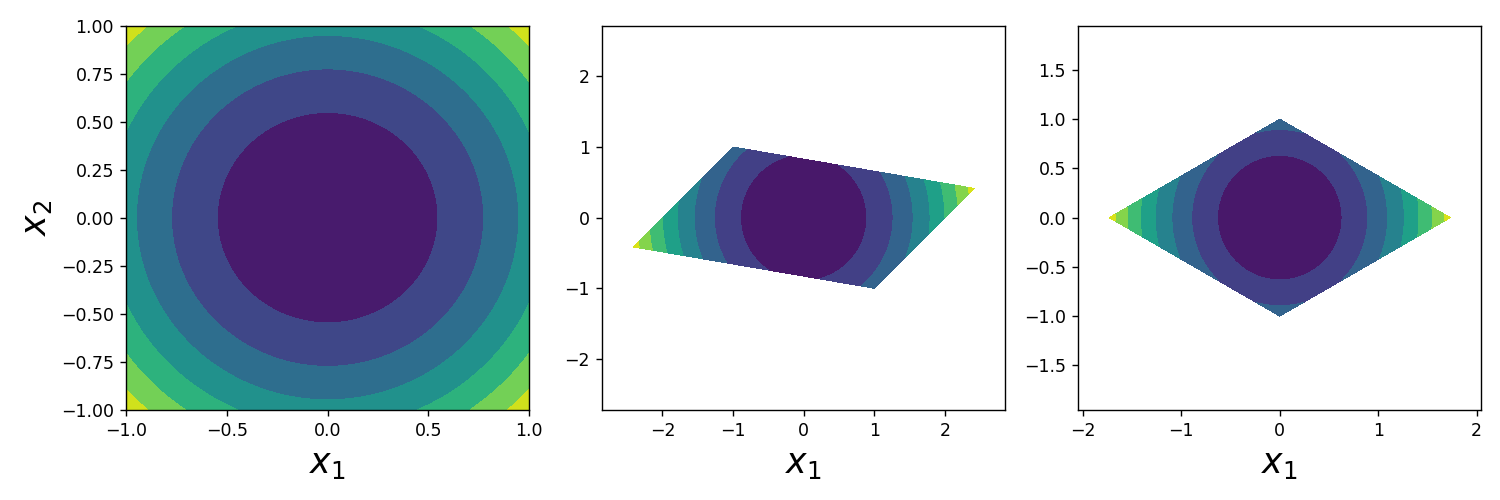}
    \caption{Contour plots of quadratic functions with Hessians given by $\boldsymbol{A}_1$, $\boldsymbol{A}_2$ and $\boldsymbol{A}_3$, after applying the proposed transformation scheme to the original domain.}
    \label{fig:QuadTrans}
\end{figure}

To start, we propose a transformation scheme when the Hessian is known. First, the eigenvectors and eigenvalues of the Hessian are computed. The domain is then rotated using the eigenvectors of the Hessian and scaled by the square root of the eigenvalues for each direction. Figure~\ref{fig:QuadTrans} shows the contours using this transformation scheme, for the 3 different $\boldsymbol{A}$ matrices.
Clearly, by taking into account the curvature in all directions the problem can be recast into a domain where the function is isotropic.

\begin{figure}[h!]
    \centering
    \includegraphics[width = \columnwidth]{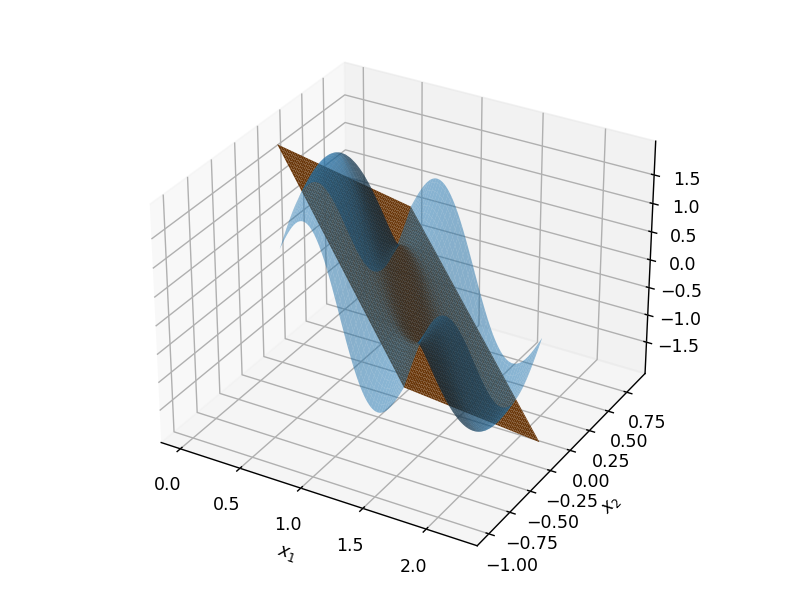}
    \caption{Global quadratic fit of the underlying function in orange and blue respectively.}
    \label{fig:GlobalFit}
\end{figure}

\begin{figure}[h!]
    \centering
    \includegraphics[width = \columnwidth]{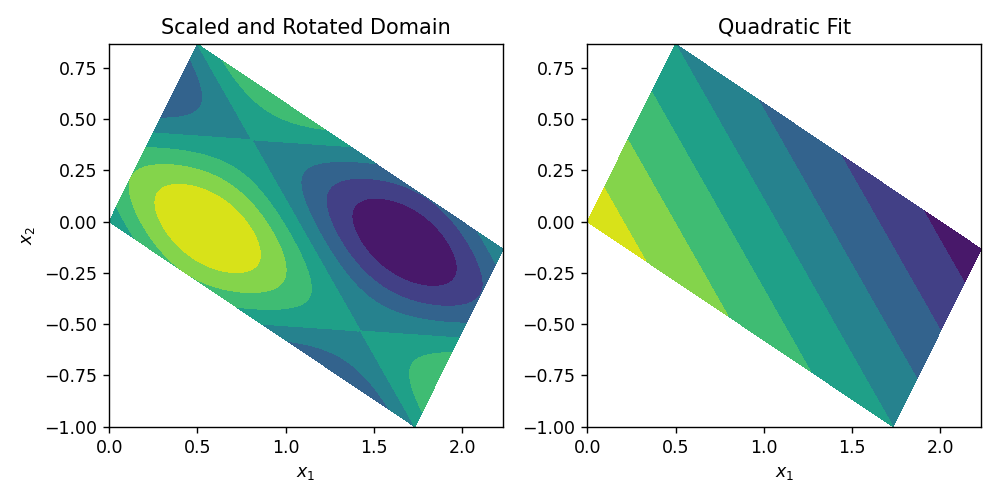}
    \caption{Contour plots of the global quadratic fit and the underlying Function.}
    \label{fig:GlobalFitContour}
\end{figure}

This scheme must now be generalised such that it can be implemented on any non-linear function (i.e.\ Hessian unknown). Initially, it seems reasonable to take some global curvature measure as, after all, the surrogate is fit on the entire domain. The issue with this assumption is demonstrated on the problem in Equation~(\ref{eq:NumericalExample}) in the scaled and rotated domain. Figure~\ref{fig:GlobalFit} shows the function overlaid with a quadratic fit of the function, while Figure~\ref{fig:GlobalFitContour} shows the contour plots of the quadratic fit and the function.

The global quadratic fit offers very little resemblance to the curvature of the underlying function. This occurs as the quadratic assumption cannot capture the full non-linearity of the underlying function across the entire domain. The regressed quadratic fit instead offers a poor representation of the underlying curvature as it completes a global least squares fit using function information. Although the regressed quadratic fit has a low function value error, as this is the information it is constructed with, it offers a poor representation of the curvature of the underlying function. This paper, therefore, proposes the use of \emph{local} quadratic fits to inform a \emph{global} curvature-based transformation scheme. 

\begin{figure}[h!]
    \centering
    \includegraphics[width = \columnwidth]{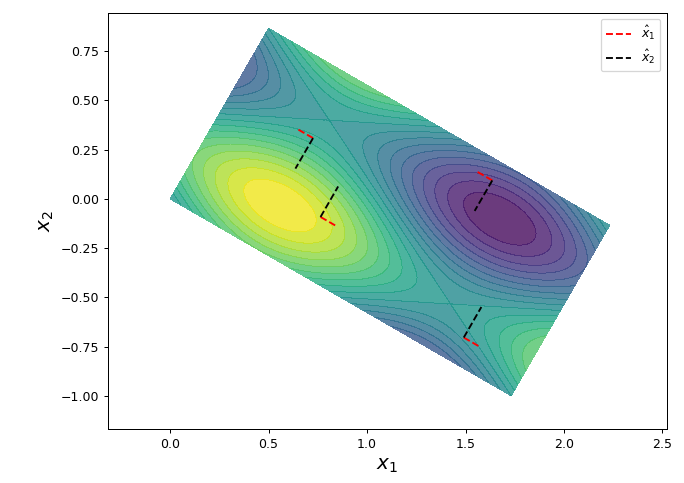}
    \caption{Local Hessian information overlaid with the rotated underlying function.}
    \label{fig:LocalHess}
\end{figure}

Figure~\ref{fig:LocalHess} shows the rotated underlying function overlaid with the eigenvalues and eigenvectors of the Hessians, i.e.\ local curvature information, at random locations in the domain. As before, the dashed lines indicate the eigenvectors, and their lengths are selected in proportion to the corresponding eigenvalues. At some locations, the eigenvalues are similar in both directions, but it occurs more frequently than the eigenvalue in the rotated $x_2$ direction (black dashed line) is larger than the eigenvalue in the rotated $x_1$ direction (red dashed line). Therefore the local curvature information more accurately reflects the curvature of the underlying function, rather than the curvature of a global approximation of the underlying function.

\begin{figure}[h!]
    \centering
    \includegraphics[width = \columnwidth]{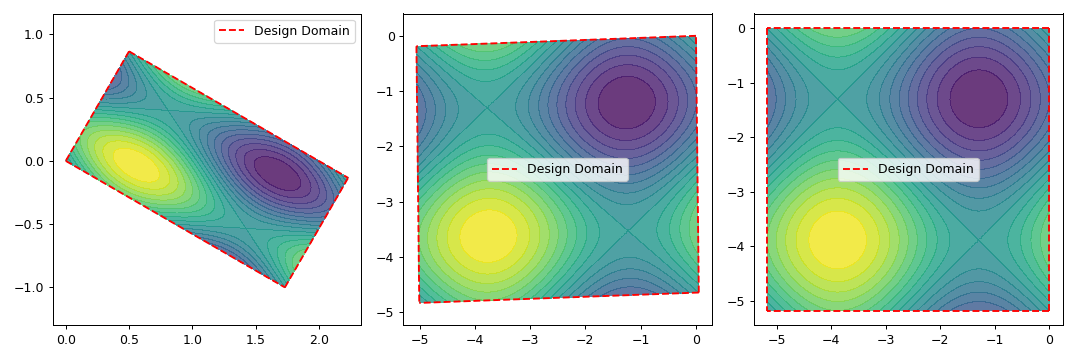}
    \caption{Contour plots of Equation~(\ref{eq:NumericalExample}) in the scaled and rotated domain, a transformed domain at 5 samples, and a transformed domain at 9 samples.} 
    \label{fig:TransSchemeRes}
\end{figure}

To remove variance in the local information some average measure of the local measures must be found. Obtaining an average orthogonal matrix from all the local eigenvectors is not a trivial computation \cite{Liski2012}. Averages of orthogonal matrices are not themselves orthogonal. Therefore  one average \emph{global} Hessian is created from the many \emph{local} Hessians. This is done by using the decomposition used in the Saddle-Free Newton method \cite{Dauphin2014}
\begin{equation}
    \boldsymbol{H} = \boldsymbol{V} \boldsymbol{\Sigma} \boldsymbol{V}^{\tiny{T}},
    \label{eq:Sadle}
\end{equation}
where $\boldsymbol{V}$ and $\boldsymbol{\Sigma}$ are the eigenvectors and a diagonal matrix containing the eigenvalues along the diagonal, respectively. Each local Hessian is then recreated by taking the absolute value of the eigenvalue matrix,
\begin{equation}
    \boldsymbol{H}_{\mbox{\tiny new}} = \boldsymbol{V} | \boldsymbol{\Sigma} | \boldsymbol{V}^{\tiny{T}}.
    \label{eq:Reconstruct}
\end{equation}
The average \emph{global} Hessian matrix is then calculated by taking the average of these reconstructed \emph{local} Hessians: 
\begin{equation}
    \boldsymbol{H}_{\mbox{\tiny avg}} = \frac{1}{N} \sum_i^N \boldsymbol{H}_{\mbox{\tiny new}}.
    \label{eq:AvgHess}
\end{equation}
Next the eigenvalues and eigenvectors of this average global Hessian is computed. The domain is rotated using the the eigenvectors as columns in an orthogonal matrix, and each direction is scaled with the square root of the eigenvalues. To demonstrate the proposed method Figure~\ref{fig:TransSchemeRes} shows contour plots of Equation~(\ref{eq:NumericalExample}) in the scaled and rotated domain, a transformed domain computed from 5 random samples, and a transformed domain computed from 9 random samples. Although in this example the local Hessians are known analytically, in general, the local Hessians must be estimated from data.

\subsection{Hessian Estimation}

What is clear from arguments presented in Sections \ref{sec:Charateristics} and \ref{sec:TransScheme} is that some understanding of the nature of the curvature of the underlying function is required. Although the exact nature of the curvature is often unavailable, there are many methods that can estimate the curvature or Hessian of a function. 

For the research completed in this paper two methods are selected depending on the information available. In the case where gradient information is available, the Symmetric Rank 1 (SR1) Hessian update method \cite{OptWilke} is used:
\begin{equation}
    \boldsymbol{H}_{k+1} = \boldsymbol{H}_{k} + 
    \frac {(\boldsymbol{y}_k - \boldsymbol{H}_k \Delta \boldsymbol{x}_k)(\boldsymbol{y}_k - \boldsymbol{H}_k \Delta \boldsymbol{x}_k)^T}{(\boldsymbol{y}_k - \boldsymbol{H}_k \Delta \boldsymbol{x}_k)^T \Delta \boldsymbol{x}_k}.
    \label{eq:SR1}
\end{equation}
The initial Hessian estimate $\boldsymbol{H}_0$ is an identity matrix and the term $\boldsymbol{y}_k$ is defined as
\begin{equation}
    \boldsymbol{y}_k = \nabla F(\boldsymbol{x}_k + \Delta \boldsymbol{x}_k) - \nabla F (\boldsymbol{x}_k).
\end{equation}
To ensure that the \emph{local} Hessian approximation is rank sufficient, $N$ SR1 updates are performed at the 
$N$ closest points surrounding the point where the Hessian is estimated. This of course requires the gradient vector at each of these $N$ points.

Otherwise, if only function information is available, a quadratic function is fitted \emph{locally} to the underlying function. This quadratic function takes the form
\begin{equation}
\boldsymbol{f} = \sum_i^N \sum_j^N W_{ij} x_i x_j + \sum_k^N W_k x_k + W_c,
\label{eq:QuadFit}
\end{equation}
where the weights $W_{ij}$, $W_k$, and $W_c$ are associated with the quadratic and coupling terms, the linear terms, and the constant term in the equation respectively. The $W_{ij}$ weights solved from this fitted function can then be re-arranged into the Hessian of the quadratic fit
\begin{equation}
    \boldsymbol{H} = 
    \begin{bmatrix}
    2W_{11} & W_{12} & \hdots & W_{1n} \\
    W_{21} & 2W_{22} & \hdots & W_{2n} \\
    \vdots & \vdots & \vdots & \vdots\\
    W_{n1} & W_{n2} & \hdots & 2W_{nn} 
    \end{bmatrix},
    \label{eq:LocalHess}
\end{equation}
where $W_{12} = W_{21}$ as the matrix is symmetric. In the implementation of this paper, an interpolating fit is constructed. This requires as many function values as there are unknown coefficients in the fit. These points are selected as the closest points surrounding the point at which the Hessian is approximated, resulting in a \emph{local} approximation of the Hessian.

\begin{figure}[h!]
    \centering
    \includegraphics[width = 0.8\columnwidth]{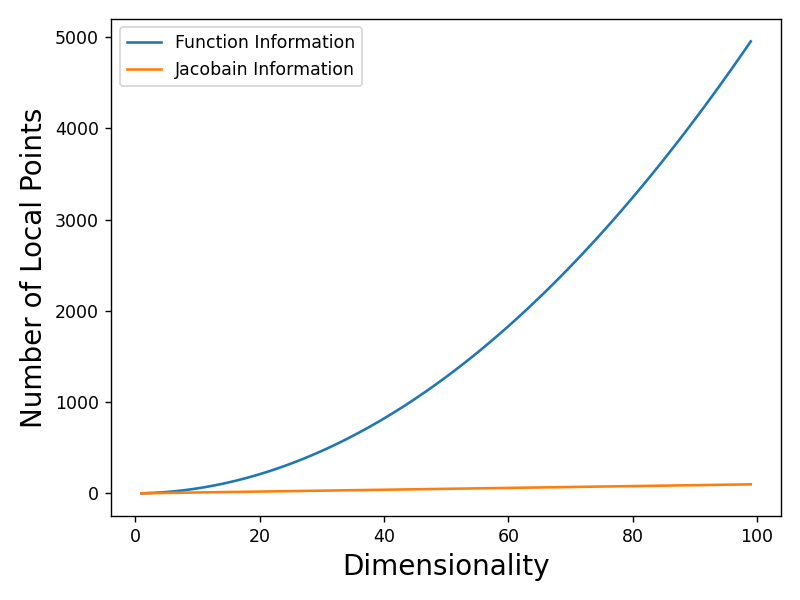}
    \caption{Number of Points required to compute a Hessian estimation as a function of the dimensionality of the problem.}
    \label{fig:LocalPoints}
\end{figure}

A key difference between these two Hessian estimation methods is the minimum number of points each method requires in order to provide an estimation of the local Hessian. The SR1 method requires $N+1$ points (the center point and the closest $N$ points) while the quadratic fit requires a local cluster containing $N(N-1)/2 + N + 1$ points in $N$-dimensional space. This implies that when gradient information is available, the proposed transformation scheme scales favorably with problem dimension (linear scaling), while the function value-based Hessian approximation method becomes prohibitively expensive (quadratic scaling). This scaling behaviour is shown visually in Figure~\ref{fig:LocalPoints}.

\subsection{Effect of Transformation on the Gradient Vector}
\label{sec:TransGrad}
When the domain is transformed, the gradients are indirectly also transformed. Therefore the gradients need to be transformed into the new domain before they are used in the construction of the surrogate in the new domain. 
This is done by first expressing the underlying function as a function of the transformed domain
\begin{equation}
    \boldsymbol{F}(\hat{\boldsymbol{x}}) = \boldsymbol{F}(\hat{\boldsymbol{x}}(\boldsymbol{x})).
    \label{Eq:Frot}
\end{equation}

Using the chain rule Equation~(\ref{Eq:Frot}) becomes 
\begin{equation}
    \frac{d \boldsymbol{F}}{d \boldsymbol{x}} = \frac{d \boldsymbol{F}}{d \hat{\boldsymbol{x}}} \frac{d\hat{\boldsymbol{x}}}{d \boldsymbol{x}}, 
\end{equation}
where $\frac{d \boldsymbol{F}}{d \boldsymbol{x}}$ is the gradients that were found when the underlying function was sampled and $\frac{d \boldsymbol{F}}{d \hat{\boldsymbol{x}}}$ is the gradients in the new transformed domain. Therefore the new gradient vector can be found in solving
\begin{equation}
    \frac{d \boldsymbol{F}}{d\hat{\boldsymbol{x}}} = \frac{d \boldsymbol{F}}{d \boldsymbol{x}} \left( \frac{d \hat{\boldsymbol{x}}}{d \boldsymbol{x}} \right)^{-1},
    \label{eq:NewGrad}
\end{equation}
The required term $(\frac{d\hat{\boldsymbol{x}}}{d \boldsymbol{x}})^{-1}$ follows from 
\begin{equation}
    \hat{\boldsymbol{x}} = \boldsymbol{R} \boldsymbol{S} \boldsymbol{x}.
    \label{Eq:Xrot}
\end{equation}
Taking the gradient of Equation~(\ref{Eq:Xrot}) yields
\begin{equation}
    \frac{d \hat{\boldsymbol{x}}}{d \boldsymbol{x}} = \boldsymbol{R} \boldsymbol{S}.
\end{equation}
Since $\boldsymbol{R}$ is a orthogonal matrix, $\boldsymbol{R}^{-1} = \boldsymbol{R}^\intercal$. Therefore, 
\begin{equation}
\left( \frac{d \hat{\boldsymbol{x}}}{d \boldsymbol{x}} \right)^{-1} = \left( \boldsymbol{R} \boldsymbol{S} \right)^{-1} = \boldsymbol{S}^{-1} \boldsymbol{R}^{-1} = \boldsymbol{S}^{-1} \boldsymbol{R}^\intercal.
\end{equation}
Since the scaling matrix $\boldsymbol{S}$ is a diagonal matrix, its inverse is simply the inverse of each diagonal entry placed in the same location on the diagonal. The final transformed gradient from Equation~(\ref{eq:NewGrad}) then becomes
\begin{equation}
    \frac{d \boldsymbol{F}}{d \hat{\boldsymbol{x}}} = \frac{d \boldsymbol{F}}{d \boldsymbol{x}} \boldsymbol{S}^{-1} \boldsymbol{R}^\intercal.
    \label{eq:FullGradTrans}  
\end{equation}

\subsection{Summary of Proposed Transformation Procedure}

The implementation of the proposed transformation procedure can be separated into 3 Sub-procedures. The first Sub-procedure, Sub-procedure \ref{al:Trans}, iterates through all the sampled points and calculates an average Hessian estimation. 

\begin{algorithm}[h!]
\SetAlgoLined

\SetKwInOut{Input}{Input}
\SetKwInOut{Output}{Output}

\Input{Sampled Information of the Underlying Function.}
\Output{Transformation Matrix Rotation $\boldsymbol{R}$, and scaling matrix $\boldsymbol{S}$.}

 \For{All sampled points}{
 \eIf{Gradient Information is available}{
 Use Procedure \ref{al:GradHess}
 }{
 Use Procedure \ref{al:FuncHess}}
 }
 \For{All Hessian Estimations}{
 Compute Equation~(\ref{eq:Reconstruct})
 }
 Compute Equation~(\ref{eq:AvgHess})\;
 \textbf{Return} The eigenvectors and eigenvalues of the average Hessian.
\caption{Transformation Procedure}
\label{al:Trans}
\end{algorithm}

\begin{figure}
    \includegraphics[width = \columnwidth]{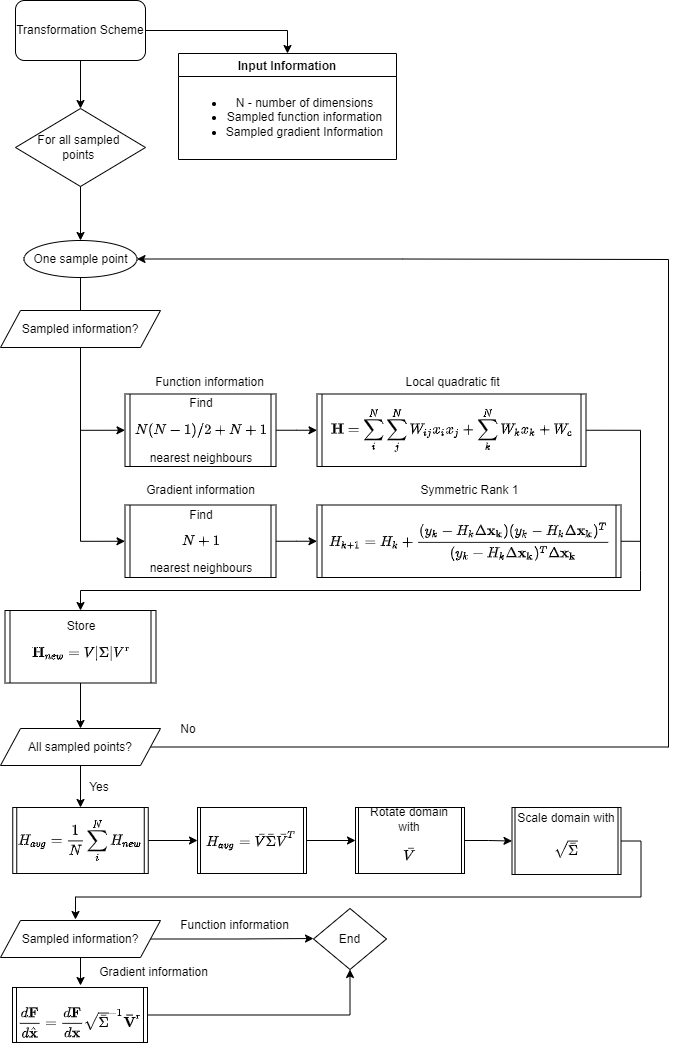}
    \caption{Flow diagram summarising the proposed transformation scheme.}
    \label{fig:Flow}
\end{figure}

Sub-procedures \ref{al:GradHess} and \ref{al:FuncHess} compute local Hessian estimations from some subset of points in the sample set. This summary is presented visually as a flow diagram in Figure~\ref{fig:Flow}.

\begin{algorithm}[h!]
\SetAlgoLined
\SetKwInOut{Input}{Input}
\SetKwInOut{Output}{Output}

\Input{A sampled point}
\Output{A Local Hessian Estimation}

 Find the $n + 1$ closest points\;
 Arrange from furthest to closest\;
 \For{Closest Points Subset}{
 Compute Equation~(\ref{eq:SR1})
 }
\textbf{Return} The local Hessian Estimation.
\caption{Gradient Information Based Hessian Estimation}
\label{al:GradHess}
\end{algorithm}

\begin{algorithm}[h!]
\SetAlgoLined
\SetKwInOut{Input}{Input}
\SetKwInOut{Output}{Output}

\Input{A sampled point}
\Output{A Local Hessian Estimation}

Find the $2n + 1$ closest points\;
Fit local Quadratic function using Equation~(\ref{eq:QuadFit}) \;
Rearrange weight vector into Hessian using Equation~(\ref{eq:LocalHess})\;
\textbf{Return} The local Hessian Estimation.
\caption{Function Information Based Hessian Estimation}
\label{al:FuncHess}
\end{algorithm}

\section{Test Problem}

In order to further evaluate i) the benefit of adequate domain transformation, and ii) the proposed transformation scheme, an $N$-dimensional test problem is constructed. The test problem will then be used to investigate the benefit of appropriate domain transformation as a function of problem dimension. If we select the test function as a decomposable function 
\begin{equation}
    f(\boldsymbol{x}) = f_1(x_1) + f_2(x_2) + \dots + f_n (x_n),
\end{equation}
then the resulting Hessian will be a diagonal matrix. Then independent scaling along each coordinate axis might create an isotropic or near-isotropic function. Therefore we select our test function as a decomposable function, ensuring that we know the optimal reference frame in which to express the function. The remaining feature that we deliberately embed into the test function, is varying length scales in different coordinate directions.
This results in a test function for which we can easily alter certain characteristics, such as problem dimension and complexity. The fact the key characteristics of the function can be easily altered allows for an independent study of desired characteristics without the need to create a new test function entirely. The test function is chosen to have the form
\begin{equation}
    \boldsymbol{f}(\boldsymbol{x}) = \frac{1}{N} \sum_{i = 1}^N A_i \sin(F_i x_i),
    \label{eq:TestFunc}
\end{equation}
where $N$ is the problem dimension and $F_i$ and $A_i$ are the frequency and amplitude in the $i^{th}$ coordinate direction. The amplitudes and frequencies are found from
\begin{align}
    A_i & = -2\exp{\frac{-(2i -N)^2}{N}} + 3, \label{eq:ConsAmp} \\
    F_i & = \frac{3\pi}{2 + 2\exp{\frac{-20i + N}{2}}} + \frac{\pi}{2}. \label{eq:ConsFreq}
\end{align}
These frequency and amplitude equations attempt to keep the complexity of the function relatively constant as the problem dimension increases. The frequency is bound between $[0.5\pi; 2\pi]$ and the amplitude between $[1, 3]$.

Another feature that is easily added to the test function, is to rotate the problem into an arbitrary reference frame. As the original test function exhibits a diagonal Hessian, a rotation of the design space is added to create a problem where the variables are not independent. This version of the test function will then assess how well the rotation aspect of the proposed transformation scheme works. The original domain is rotated using a random rotation matrix $\boldsymbol{R}$ created from 
\begin{equation}
   \boldsymbol{R} = \mbox{expm}(\pi( \boldsymbol{A} - \boldsymbol{A}^\intercal ) ),
   \label{eq:Rotation}
\end{equation}
where $\boldsymbol{A}$ is a random matrix with elements sampled between $[-0.5, 0.5]$ and $\mbox{expm}$ is the exponential map. The exponential map of a skew matrix $(\boldsymbol{A}-\boldsymbol{A}^\intercal)$ results in an orthogonal matrix \cite{OptWilke}. This is done as during the testing phase the ideal transformed domain is available (it is assumed to be the reference frame in which the test function is a decomposable function) by simply using the same method shown in Section \ref{sec:NumericalExample}. 

In this research, the case where gradient information is available is also discussed. Therefore, the gradients of the $N$-dimensional test function are needed. The gradient of Equation~(\ref{eq:TestFunc}) is simply
\begin{equation}
    \dfrac{\partial F_i}{\partial x_i}  = A_i F_i \cos (F_i x_i),
\end{equation}
where in the case of domain rotation, the process detailed in Section \ref{sec:TransGrad} is followed.

\section{Results}
\label{sec:Results}
The numerical results in this section follow a two-step process
\begin{enumerate}
    \item a domain transformation,
    \item followed by surrogate construction.
\end{enumerate}
The results, therefore, attempt to separate the contribution of these two steps to the performance of a surrogate. Specifically, the information used to perform domain transformation is deliberately separated from the information used to construct the surrogate. 

This is done by constructing the two types of surrogates discussed in Section \ref{sec:Construction}, the function-value based and GE surrogates, in five different domains. These five domain transformations are
\begin{itemize}
    \item Domain transformation (rotation and scaling) performed using gradient information,
    \item Domain transformation (rotation and scaling) performed using function information,
    \item Domain scaling (no rotation) using the Kriging hyper-parameter optimization strategy discussed in Section \ref{sec:HyperParameter},
    \item The domain is min-max scaled (no rotation) to $[0; 1]$ in all dimensions, and
    \item The ideal transformation is used, as discussed in Section~\ref{sec:NumericalExample}. This transformation is only possible since we have an analytical expression for the underlying function, hence we can compute the Hessian analytically.
\end{itemize}
By using two different models in five different domains, it will become apparent in the results if the domain consistently impacts the performance of the surrogate model, regardless of the information used in the construction of the model. The two surrogate models, function value and GE, have the same model flexibility, i.e.\ the same number of centres, to further isolate the effect the domain the models are constructed in has on the performance of the surrogate model. By fixing the flexibility of the surrogate model it will be shown that the ill-suitably of the domain the model is constructed in, and not a lack of construction information is the main source of the approximation error. 

As in Section \ref{sec:NumericalExample}, the surrogates are constructed using various numbers of sampled points. This construction is repeated 50 times for a fixed number of samples, allowing the calculation of the mean and variance of the surrogate performance. To offer a more visual demonstration of the effect of the domain transformation, 1D lines through $N-$dimensional space are constructed, on which both the underlying function and the surrogate are sampled. This allows simple visualisation of the higher dimensional problems.  

\subsection{RMSE Results}
\label{sec:RMSE}

\begin{figure}[h!tbp]
    \centering
    \includegraphics[width = \columnwidth]{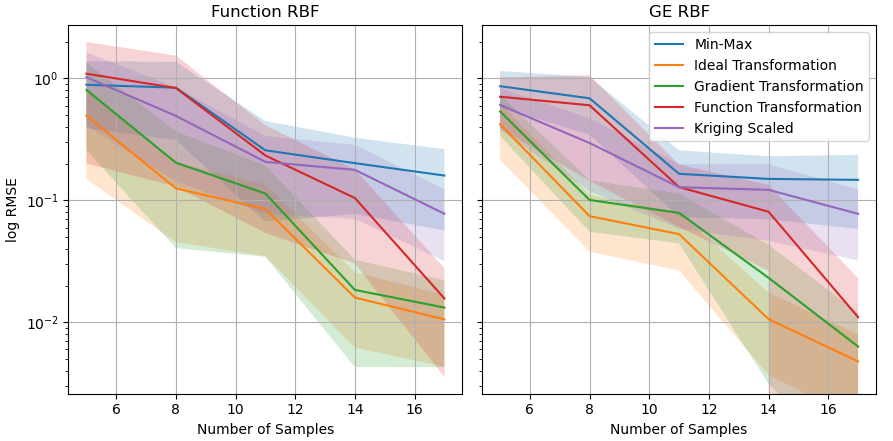}
    \caption{Log RMSE results for the 2-dimensional test problem.}
    \label{fig:2D_Results}
\end{figure}

The RMSE of the surrogates is found by sampling the error at $10^5$ test points. Such a large number of test points is selected to ensure that an accurate RMSE is computed even for the high-dimensional versions of the test problem. This process is then repeated 50 times to be able to compute the average RMSE error, as well as the variance in the RMSE. Figure \ref{fig:2D_Results} presents the results for the 2-dimensional test problem. The average RMSE (solid lines) and the variance in RMSE (shaded areas) are shown for the functional and GE RBFs in all five construction domains. The RMSE results are presented in the log domain so that the performance of the models can be compared across a wide range of accuracy levels. 

This 2D example shows that there is a benefit in constructing the surrogate in the transformed domain instead of the $[0;1]$ scaled domain, most noticeable when only 7 or 8 samples are used for the surrogate construction. It is also noticeable that only scaling the domain, i.e.\ the Kriging scaled results, is not nearly as beneficial as complete domain transformation (scaling and rotation). 

\begin{figure}[h!tbp]
    \centering
    \includegraphics[width = \columnwidth]{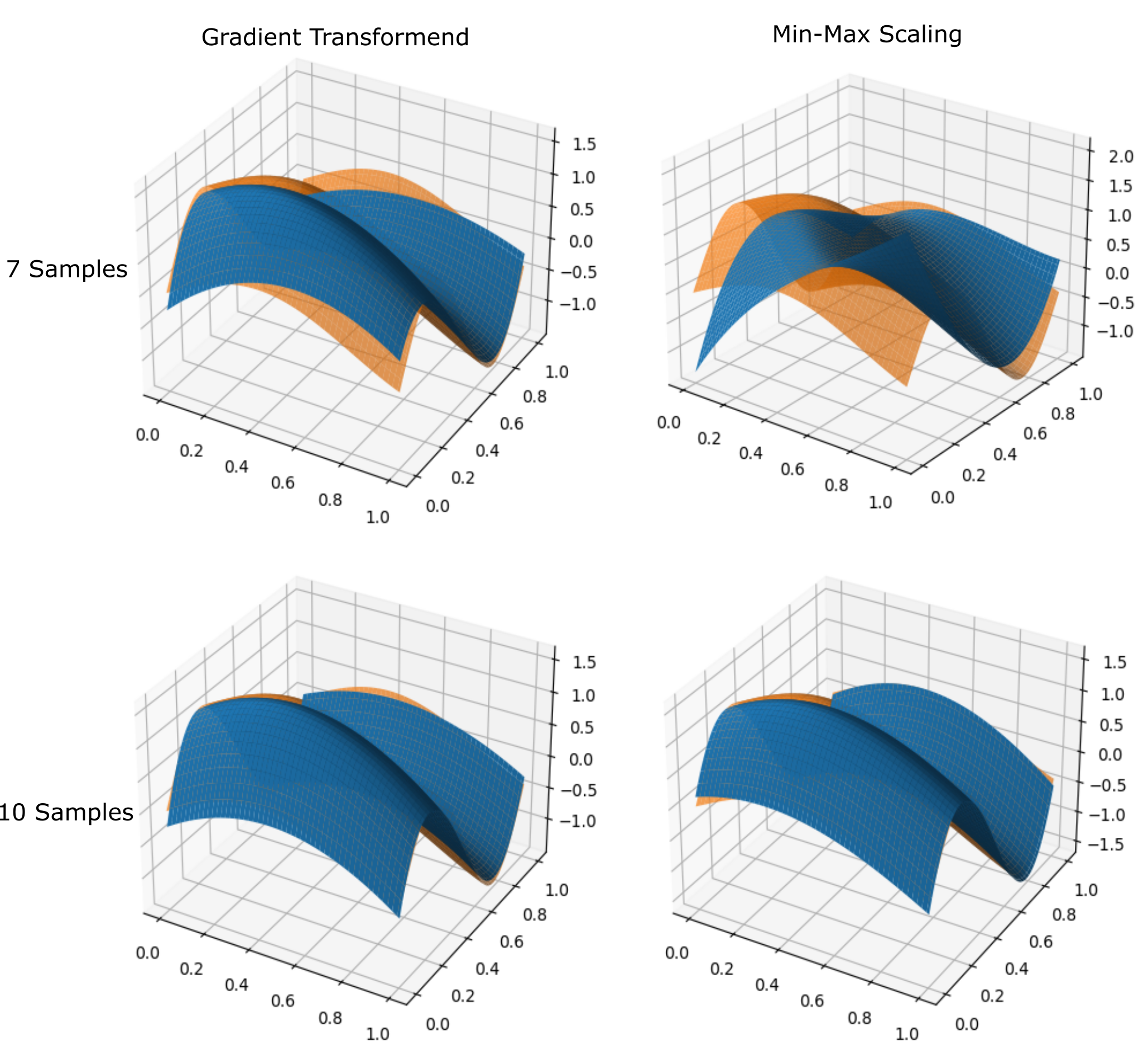}
    \caption{3D Plots showing the benefit of appropriate domain transformation with the underlying function in orange and the GE surrogate models in blue.}
    \label{fig:2D_Examples}
\end{figure}

As soon as a sufficient number of samples is used (between 10 and 12 in this test problem), the sampling density is sufficient to overcome the non-isotropic and coupled nature of the problem and constructing the surrogate in any of the five domains returns satisfactory results. This is demonstrated in  Figure~\ref{fig:2D_Examples} using the 2D GE surrogate models for 7 or 10 samples. 

\begin{figure}[h!tbp]
    \centering
    \includegraphics[width = \columnwidth]{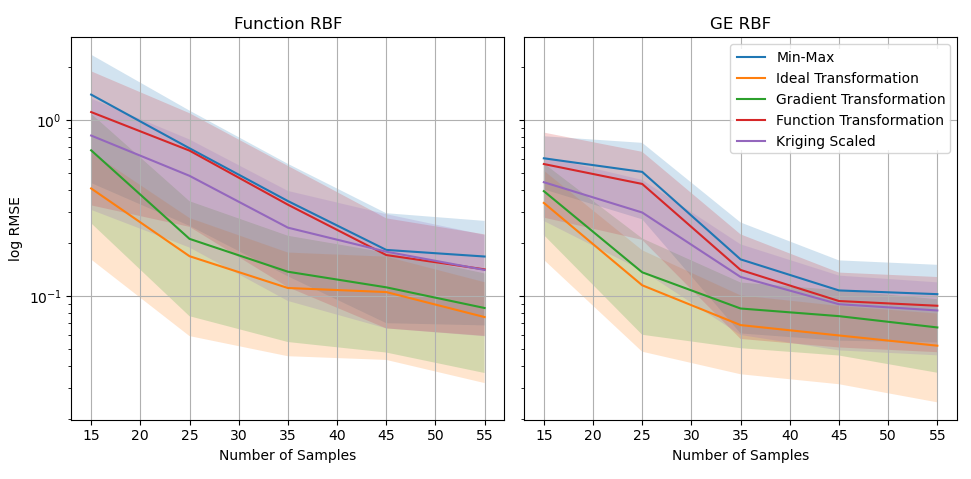}
    \caption{Log RMSE results for the 4-dimensional test problem.}
    \label{fig:4D_Results}
\end{figure}

\begin{figure}[h!tbp]
    \centering
    \includegraphics[width = \columnwidth]{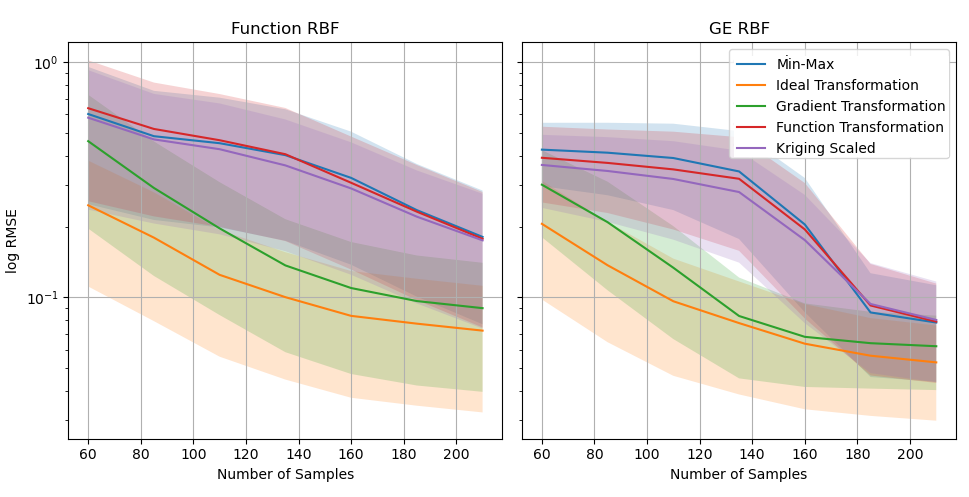}
    \caption{Log RMSE results for the 8-dimensional test problem.}
    \label{fig:8D_Results}
\end{figure}

As will be shown, overcoming the coupled and non-isotropic nature of the function with dense enough sampling becomes far more difficult in higher dimensional problems. Figures~\ref{fig:4D_Results} and \ref{fig:8D_Results} present the results for the 4 and 8-dimensional problems respectively. 

Another important aspect of the results is the rate of improvement of the surrogate. That is quantifying the improvement when additional information or samples are added to the surrogate. When Figure \ref{fig:8D_Results} is considered it is clear that initially, the improvement is minimal when no domain transformation is completed. In the case where the ideal transformation or the proposed gradient transformation scheme is implemented far more performance is gained at low sampling density. 

This increase in problem dimension highlights both the importance of a complete transformation scheme as well as the benefit of gradient information. Firstly, for the 4-dimensional problem, there is some benefit of the Kriging-based scheme over the proposed function transformation and the simple min-max scaling. But, as the problem dimension increases to 8, this benefit diminishes to almost zero. The second observation to note is the clear performance gain when a suitable completely transformed construction domain is used. This gain is evident in both the ideal domain and the gradient-transformed domain cases. Gradient information offers a better approximation of local curvature, and therefore, returns a near-optimal approximation of the ideal transformed construction domain.   

\begin{figure}[h!tbp]
    \centering
    \includegraphics[width = \columnwidth]{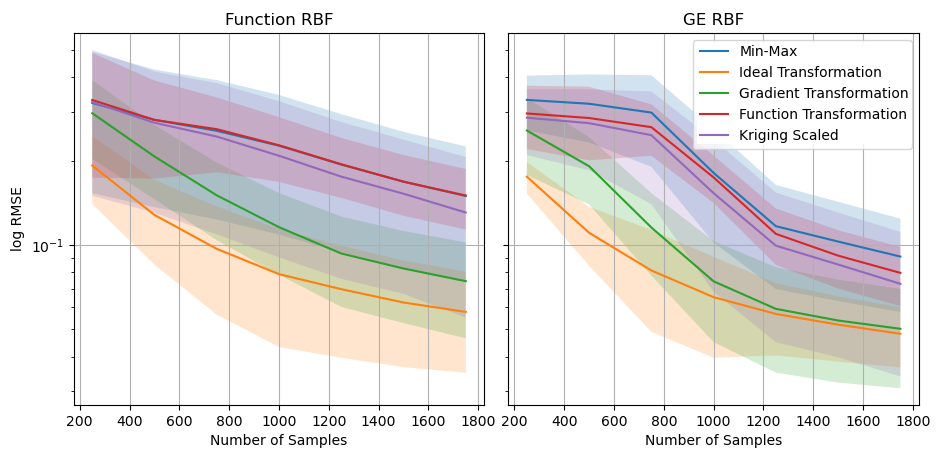}
    \caption{Results for the 16-dimensional test problem}
    \label{fig:16D_Results}
\end{figure}

The problem dimension is then further increased to 16 and the same results are repeated in Figure \ref{fig:16D_Results}. From these results, it becomes apparent that the benefit of appropriate complete domain transformation, over both min-max scaling or Kriging scaling, grows with problem dimension. As with the lower dimensional problems, the surrogates constructed in ill-suited domains offer minimal performance improvement in low sample density scenarios when additional samples are added. This slow rate of improvement for the ``non-transformed'' surrogate means that the proposed gradient-based transformation scheme and the ideal transformation domain require far less computational cost to achieve the same accuracy. For example, if the 16-dimensional problem had a goal RMSE of 0.1 the proposed transformation scheme would require, on average, 1200 and 800 samples for the function and GE models respectively, while the standard min-max scaling would require $>$2000 and 1600 samples. Therefore, for this simple test function, the proposed transformation scheme results in almost half the computational cost of the standard scaling procedure. 

\subsection{1D lines in $N-$dimensional space}

\begin{figure}[h!tbp]
    \centering
    \includegraphics[width = \columnwidth]{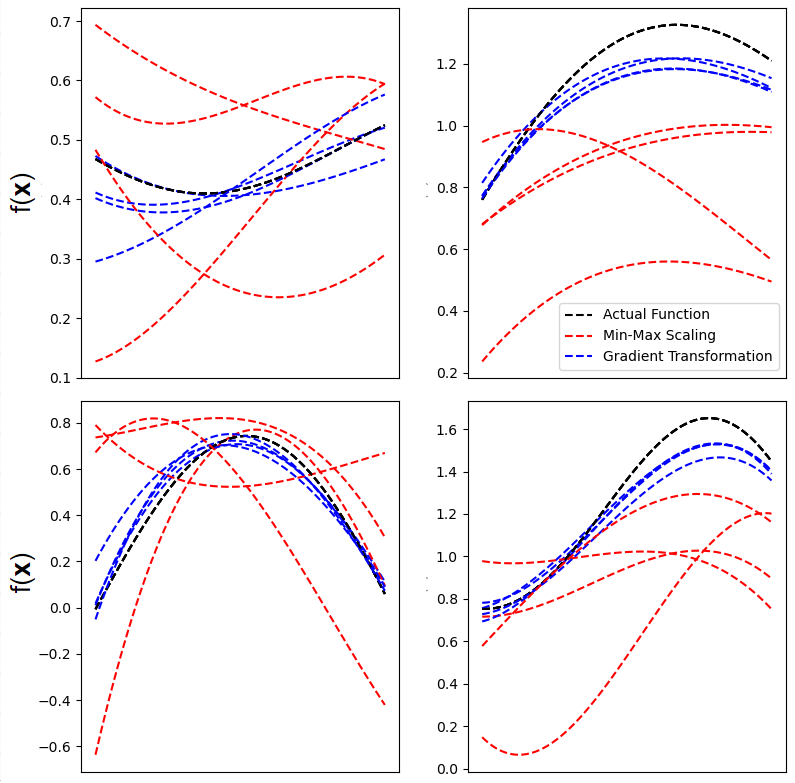}
    \caption{1D lines in 4D space: 30 Samples used to construct GE-RBFs.}
    \label{fig:4D_LowSample}
\end{figure}

\begin{figure}[h!tbp]
    \centering
    \includegraphics[width = \columnwidth]{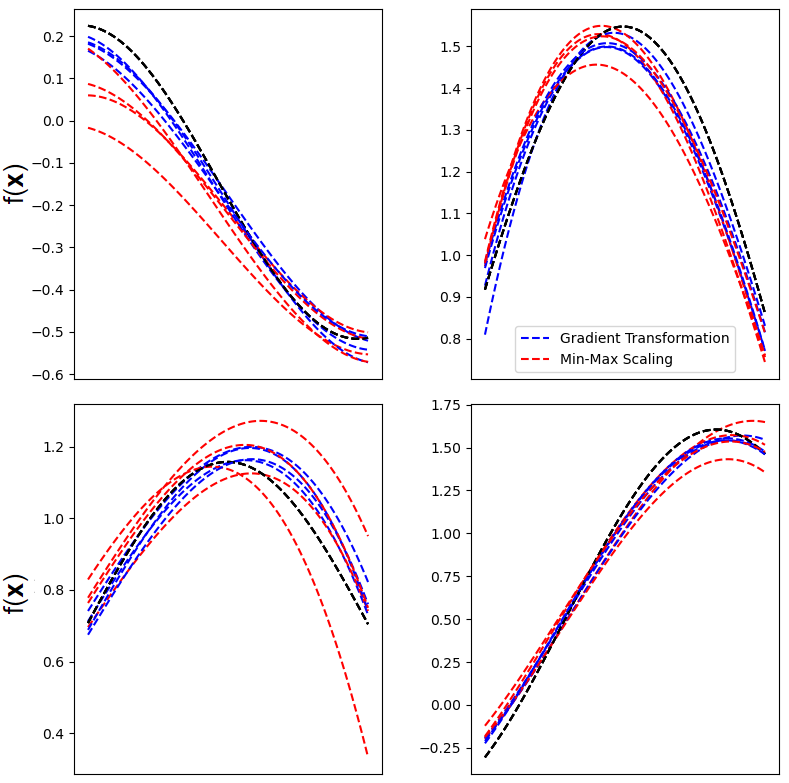}
    \caption{1D lines in 4D space: 50 Samples used to construct GE-RBFs.}
    \label{fig:4D_HighSample}
\end{figure}

For each test problem dimension, four 1D lines through $N-$dimensional space are sampled for GE-RBFs in the min-max scaled domain, and in the gradient-transformed domains. The proposed gradient transformation scheme is compared to the min-max scaling case since this is the {\it de facto} standard in RBF construction. For each domain, four lines are sampled simply to demonstrate the possible variance in the shapes of the surrogates. The number of samples are then increased, such that each test problem dimension will have two sets of figures, to visualise the improvement of the surrogates as more samples are added. The four-dimensional case is presented in Figures \ref{fig:4D_LowSample} and \ref{fig:4D_HighSample} for 30 and 50 samples respectively.

The 1D lines through $N$-dimensional space demonstrate the difficulty the surrogates have when presented with a problem constructed in an unsuitable domain. For the standard min-max scaling surrogate approach, the curves at low sampling density offer almost no resemblance to the shape of the underlying function, and each generated surrogate can return wildly different results. 

At low sample density, the min-max scaled surrogates struggle to offer any resemblance to the shape of the underlying function and are therefore greatly dependent on the randomness in the sampling locations. When the number of samples is increased the performance of both surrogates improves, but, as with the RMSE results, the transformed surrogates remain more accurate and consistent. 

\begin{figure}[h!tbp]
    \centering
    \includegraphics[width = \columnwidth]{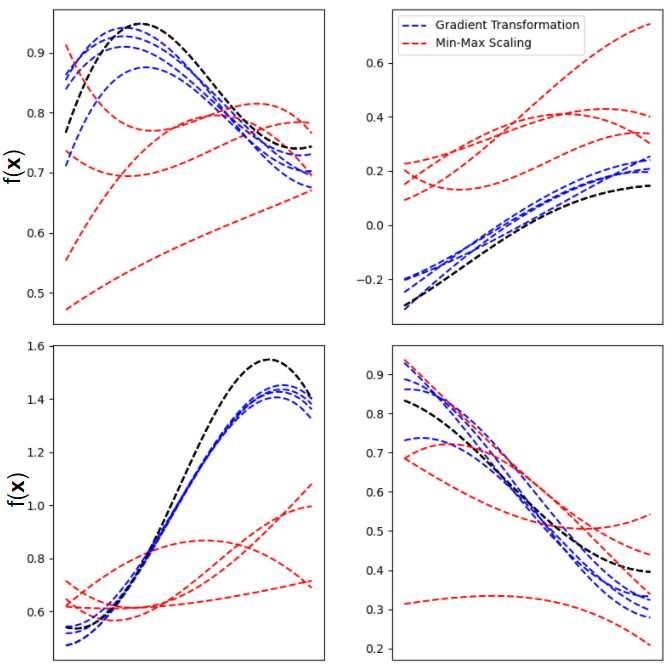}
    \caption{8-dimensional lines at 130 Samples for GE-RBFs.}
    \label{fig:8D_LowSample}
\end{figure}

\begin{figure}[h!tbp]
    \centering
    \includegraphics[width = \columnwidth]{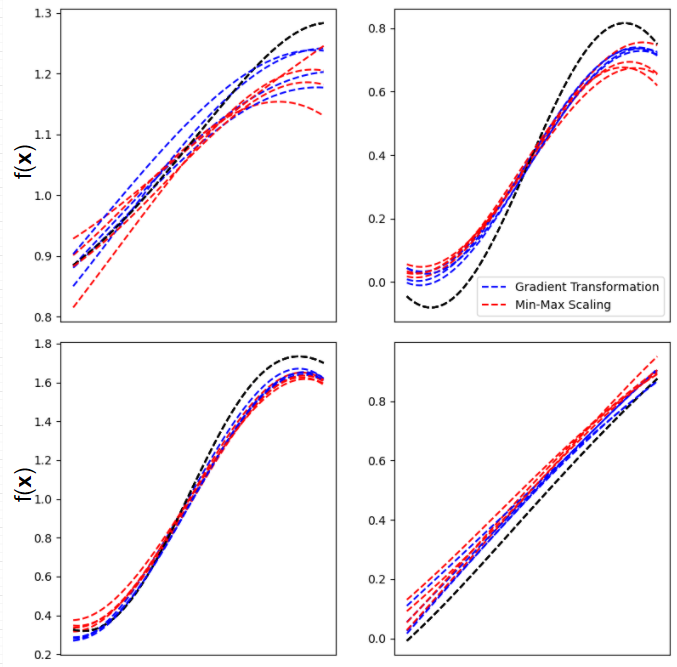}
    \caption{8-dimensional lines at 190 Samples for GE-RBFs.}
    \label{fig:8D_HighSample}
\end{figure}

\begin{figure}[h!tbp]
    \centering
    \includegraphics[width = \columnwidth]{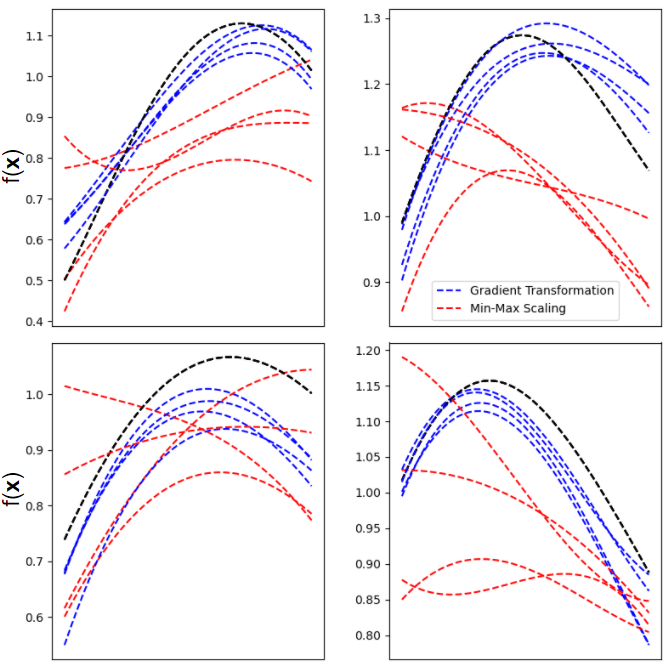}  
    \caption{16-dimensional lines at 850 Samples for GE-RBFs.}
    \label{fig:16D_LowSample}
\end{figure}

\begin{figure}[h!tbp]
    \centering
    \includegraphics[width = \columnwidth]{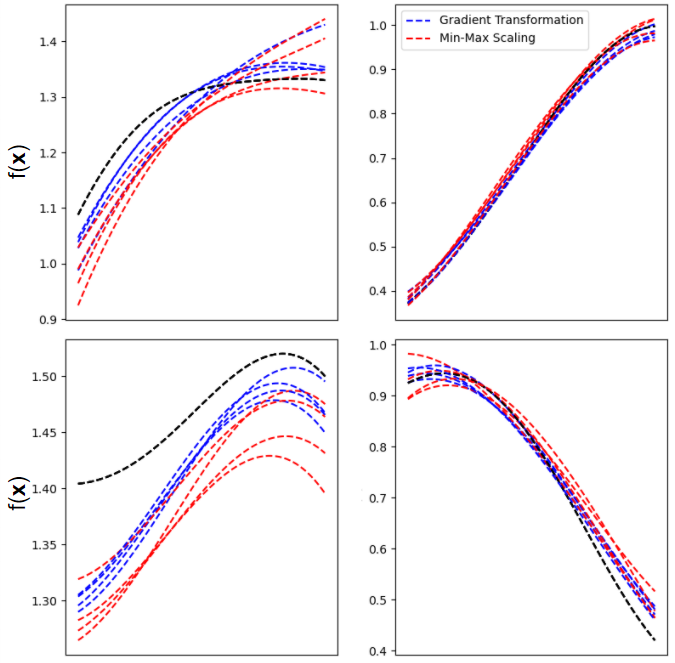}
    \caption{16-dimensional lines at 1250 Samples for GE-RBFs.}
    \label{fig:16D_HighSample}
\end{figure}

The testing procedure is then repeated for the 8 and 16-dimensional cases in Figures \ref{fig:8D_LowSample} - \ref{fig:16D_HighSample} respectively.
The same behaviour that is found in the 4-dimensional case is present in the 8 and 16-dimensions problems. At low sampling densities, if the surrogate is not constructed in an appropriate domain, the results are inconsistent and poor when compared to a surrogate constructed in a suitably transformed domain. 

\section{Conclusion}

The work presented in this paper demonstrates that the domain in which common surrogate models are constructed can have a significant influence on the predictive performance of the surrogate. This is done with a few main findings. 

Firstly, the addition of gradient information into the construction of a surrogate model will not result in the expected improvement of the predictive performance of the surrogate if the domain is not suitable. Therefore, attention needs to be given to a pre-processing step that will adequately transform the domain in which the surrogate model will be constructed. 

The information needed to inform the pre-processing step is a collection of local curvature information rather than one global estimation of the curvature. This local curvature will need to be estimated in most practical engineering problems. This estimation can be done with either gradient or function information but gradient information offers a more efficient and accurate approximation of the local curvature.

The domain the models are constructed in impacts the performance of the surrogate model regardless of the information used to construct the model. There is improvement in both the function-value and GE surrogate models when the domain the models are constructed in is transformed using the developed domain transformation scheme. The transformation must be a fully coupled rotation and scaling as only scaling the domain is not sufficient. 

Lastly, the use of gradient information allows for the estimation of local curvature to complete a powerful, automatic, and fully coupled domain transformation scheme that results in near-optimal performance. Therefore, using the gradient information to transform the domain can be far more beneficial to surrogate performance than including this information directly in the construction of the surrogate model. 

\section{Future Work}

Although the proposed transformation scheme offers a significant improvement over the standard min-max scaling scheme, there are two main scenarios that were not investigated: 
\begin{itemize}
    \item The underlying function curvature varies greatly along a principal direction (commonly referred to as non-stationary problems), and
    \item one dimension is sampled more densely than the other dimensions, such as with time series data.
\end{itemize}
These scenarios may require adaptation to the proposed transformation scheme to achieve the same level of improvement as demonstrated in this paper.

\section*{Conflict of interest}
The authors declare that they have no conflict of interest.

\section*{Replication of results}

All necessary algorithms and problem parameters for possible replication of all result presented in this work have been detailed and referenced.

\section*{Funding Sources}

There are no funding sources for this research.

\bibliographystyle{IEEEtran}      
\bibliography{Bib}   

\end{document}